\documentstyle[12pt]{article}
\begin{document}
\setlength{\textheight}{574pt}
\setlength{\textwidth}{432pt}
\setlength{\oddsidemargin}{18pt}
\setlength{\topmargin}{14pt}
\setlength{\evensidemargin}{18pt}
\newtheorem{theorem}{Theorem}[section]
\newtheorem{lemma}{Lemma}[section]
\newtheorem{corollary}{Corollary}[section]
\newtheorem{conjecture}{Conjecture}
\newtheorem{remark}{Remark}[section]
\newtheorem{definition}{Definition}[section]
\newtheorem{problem}{Problem}
\newtheorem{example}{Example}[section]
\newtheorem{proposition}{Proposition}[section]
\title{{\bf NUMERICALLY TRIVIAL FIBRATIONS}}
\date{October 24, 2000}
\author{Hajime Tsuji}
\maketitle
\begin{abstract}
We develop an intersection theory for  pseudoeffective 
singular hermitian line bundles (cf. Definition 2.4)
 on a smooth projective variety and irreducible curves on the variety.  
And we prove the existence 
of a   natural rational fibration structure associated with 
the singular hermitian line bundle. 
Also for any pseudoeffective line bundle  on a smooth 
projective variety, we prove the existence of a  natural rational 
fibration structure associated with the line bundle. 

We also characterize a numerically trivial singular hermitian 
line bundle on a smooth projective variety.  MSC32J25  
\end{abstract}
\tableofcontents
\section{Introduction}

Let $X$ be a smooth projective variety and let 
$L$ be a line bundle on $X$. 
It is fundamental to study the ring 
\[
R(X,L) := \oplus_{m\geq 0}\Gamma (X,{\cal O}_{X}(mL))
\]
(in more geometric language to 
study the Iitaka fibration associated with $L$) in algebraic geometry.
In most case, to show {\em the nonvanishing} ,i.e., \\ $\Gamma (X,{\cal O}_{X}(mL))\neq 0$ for some $m >0$ is a central problem.  

Because $R(X,L) \simeq$ {\bf C}, if $L$ is not pseudoeffective (cf. Definition 2.4),
the problem is meaningful only when  $L$ is pseudoeffective.

If $L$ is big, then for a sufficiently large $m$, 
the linear system $\mid mL\mid$ gives a  birational 
rational embedding of $X$ into a projective space. 
But if $L$ is not big, there are very few tools 
to study $R(X,L)$ except Shokurov's nonvanishing theorem \cite{sh}. 
Moreover even if $L$ is big, to study $R(X,L)$ we 
often need to study the restriction of $R(X,L)$ 
on the  subvarieties on which the restriction of 
$L$ is not big (e.g. \cite{tu3}). 

When $L$ is not big, a natural approach is to distinguish 
the  {\em  null direction} of $L$. 
Then we may consider that $L$ has positivity in 
the transverse direction. 

If $L$ has a $C^{\infty}$-hermitian metric $h$ 
such that the cuvature form $\Theta_{h}$ is semipositive, 
the {\em null foliation}
\[
\cup_{x\in X}\{ v\in TX_{x} \mid \Theta_{h}(v,\bar{v}) = 0\}
\]
defines a $C^{\infty}$-foliation on the open subset 
where the rank of the semipositive form  $\Theta_{h}$ is maximal  
and every leaf is a complex submanifold on the set. 
In this case the null direction is given by this foliation. 

But in general, a  pseudoeffective line bundle on a smooth projective variety 
does not admit a $C^{\infty}$-hermitian metric 
with semipositive curvature, even if it is nef,
although it admits a singular  hermitian metric
with positive curvature current\footnote{Here  we note that ``positive"  does not mean strict positivity (cf. Definition 2.2).  This terminology may be misleading for algebraic geometers.  For this reason I include a subsection which summarize the notion of closed positive currents. }.
Hence we need to consider a singular hermitian metric on $L$ 
in order to study $R(X,L)$. 
 
In this paper we develop an intersection  theory for singular hemitian line 
bundles with positive curvature current and curves on a smooth projective variety. 
The new intersection number measures 
the intersection of the {\em positive part} of the singular hermitian 
line bundle and the curve. 
This intersection theory is not cohomological.
 
We obtain a natural rational fibration structure in terms of 
this intersection theory as follows. 
\begin{theorem}(Fibration theorem)
Let $(L,h)$ be a pseudoeffective singular hermitian line bundle 
(cf. Definition 2.4) on a smooth projective variety $X$. 
Then there exists a unique (up to birational equivalence) rational fibration 
\[
f : X -\cdots \rightarrow Y
\]
such that 
\begin{enumerate}
\item $f$ is regular over the generic point of  $Y$, 
\item for every very general fiber $F$, 
 $(L,h)\mid_{F}$ is well defined and is numerically trivial (cf. Definition 2.9,2.10),
\item $\dim Y$ is minimal among such fibrations.
\end{enumerate}
\end{theorem}
We call the above fibration $f : X -\cdots\rightarrow Y$ 
the {\bf numerically trivial fibration} associated with $(L,h)$.
\begin{remark}
Let $X$,$(L,h)$ be as above. Then for any smooth divisor $D$ 
on $X$,
there exists a numerically trivial fibration 
\[
f_{D} : D -\cdots \rightarrow W.
\]
This is simply because the restriction of the intersection theory 
on $D$ exists (cf. Section 2.5) and the proof of the above theorem essentially
does not require the existence of the restriction of $\Theta_{h}$
to $D$. 
\end{remark} 
\begin{remark}
By the proof of Theorem 1.1 below, we see that 
the 3-rd condition in Theorem 1.1 is equivalent to : \\
$3^{\prime}.\,$ for a very general point $x\in X$ and any irreducible  
 horizontal  curve (with respect to $f$) $C$ containing $x$, 
$(L,h)\cdot C > 0$ holds (cf. Definition 2.9 for the definition of 
$(L,h)\cdot C$). 
\end{remark}
Theorem 1.1 singles out the null direction of $(L,h)$ as fibers. 
But this direction is only a part of the null direction as is shown by the following example. 
This example also shows that in general $\dim Y$ may be strictly larger 
than the numerical dimension of $L$.
\begin{example}
Let $X$ be an irreducible quotient of the open unit bidisk $\Delta^{2}$ 
in $\mbox{\bf C}^{2}$, i.e.,
\[
X = \Delta^{2}/\Gamma ,
\]
where $\Gamma$ is an irreducible cocompact torsion free lattice. 
Let $(L,h)$ denotes the hermitian line bundle whose curvature form 
comes from the Poincar\'{e} metric on the first factor. 
Then one see that $L$ is nef and $L^{2} = 0$ holds.
In particular $L$ is not big.
In this case the null foliation of $\Theta_{h}$ is 
nothing but the projection of the fibers of the first projection and 
every leaf of the  foliation is Zariski dense (actually 
even topologically dense in usual topology) in $X$.
This implies that $L$ (and hence also $(L,h)$) is numerically positive
and the numerically trivial fibration is the identity.  
\end{example}
By using an AZD (cf. Definition 2.8, Theorem 2.4 and Proposition 2.1 below), we have the following corollary.
\begin{corollary}
Let $L$ be a pseudoeffective line bundle on 
a smooth projective variety $X$ and let $h$ be 
a canonical AZD of $L$ (cf. Section 2.3).
Then there exists a unique rational fibration 
(up to birational equivalence): 
\[
f : X -\cdots \rightarrow Y
\]
such that 
\begin{enumerate}
\item $f$ is regular over the generic point of $Y$, 
\item for every very general fiber $F$,  
 $(L,h)$ is numerically trivial on $F$.
\item $\dim Y$ is minimal among such fibrations.
\end{enumerate}
Also $f$ does not depend on the choice of the canonical AZD $h$ (see 
Proposition 2.1).
\end{corollary}
We call the above fibration $f : X-\cdots\rightarrow Y$ 
the {\bf numerically trivial fibration} associated with $L$. 

The poof of Theorem 1.1 is done by finding a dominating 
family of maximal dimensional subvarieties 
on which the restriction of $(L,h)$ is numerically trivial.
{\em The heart of the proof is to prove that this family actually  
gives a rational fibration by showing that the generic point of 
a general member of the family does not intersect other members.}

The structure of numerically trivial singular hermitian line bundles 
with positive curvature current is given  as follows.
\begin{theorem}
Let $(L,h)$ be a singular hermitian line bundle on a 
smooth projective variety $X$. 
Suppose that $\Theta_{h}$ is closed positive 
and $(L,h)$ is numerically trivial on $X$. 
Then there exist at most countably many prime divisors $\{ D_{i}\}$
and nonnegative numbers $\{ a_{i}\}$ such that 
\[
\Theta_{h} = 2\pi\sum a_{i}D_{i}
\]
holds, where we have identified each $D_{i}$ with a closed 
positive current.  
More generally  let $Y$  be a subvariety of $X$ such that the restriction  
$h\mid_{Y}$ is well defined.
Suppose that  $(L,h)$ is numerically trivial 
on $Y$. 
Then the restriction $\Theta_{h}\mid_{Y}$ is a sum of at most countably many 
prime divisors with nonnegative coefficients on $Y$. 
\end{theorem}
\begin{remark}
For a divisor $D$, the current associated 
with $D$ is often denoted by $[D]$. 
But this notation is confusing with the round down of $D$ 
in algebraic geometry. 
Hence we do not use this notation in this paper. 
\end{remark}
This paper is a byproduct of the proof of the  nonvanishing theorem 
(\cite[Theorem 5.1]{tu3}). 

In this paper, I cannot refer to applications of the above 
theorems because of the length. 
These will be published separately. 

In this paper ``{\em very general}'' means outside of at most countably many union of 
proper Zariski closed subsets and ``{\em general}'' means in the sense of usual Zariski topology. 

I intended the paper to be readable for algebraic geometers 
who are not familiar with complex analytic background. 

I would like to express hearty thanks to the referee for his 
careful reading and a lot of useful comments. 

\section{Intersection theory for singular hermitian line bundles}
In this section we define an intersection number for 
a singular hermitian line bundle with positive curvature current 
on a smooth projective variety and an irreducible curve 
on it.  
This intersection number is different from the usual intersection 
number of the underlying line bundle  and the curve. 
\subsection{Closed positive currents}
In this subsection we shall review the definition and basic
notions of closed positive $(p,p)$-currents on 
a complex manifold. 
For the general facts about the theory of currents, 
see for example \cite[Chapter 3]{g-h}.
Let $M$ be a complex manifold of dimension $n$ and let 
$A^{p,q}_{c}(M)$ denote the space of $C^{\infty}$ 
$(p,q)$-forms with compact support. 
We define a topology on $A^{p,q}_{c}(M)$ such that 
a sequence $\{\varphi_{i}\}_{i=1}^{\infty}$ in $A^{p,q}_{c}(M)$ converges, 
if and only if there exists a compact subset $K$ of 
$M$ such that $\mbox{Supp}\,\varphi_{i} \subseteq K$ holds
for every $i$ and 
$\{ \varphi_{i}\}_{i=1}^{\infty}$ converges in $C^{k}$-topology on $K$ 
for every $k$ to a $C^{\infty}$ $(p,q)$-form $\varphi_{\infty}$. 
\begin{definition}
Let $M$ be a complex manifold of dimension $n$. 
The space of $(p,q)$-currents ${\cal D}^{p,q}(M)$ on $M$ is the dual space of 
$A^{n-p,n-q}_{c}(M)$. 
We define 
\[
\partial : {\cal D}^{p,q}(M)\longrightarrow {\cal D}^{p+1,q}(M)
\]
and 
\[
\bar{\partial} : {\cal D}^{p,q}(M)\longrightarrow 
{\cal D}^{p,q+1}(M)
\]
by 
\[
\partial T(\varphi ) := (-1)^{p+q+1}T(\partial \varphi )\,\,\, (T\in {\cal D}^{p,q}(M),\varphi \in A^{n-p,n-q}(M))
\]
and
\[
\bar{\partial}T(\varphi ) := (-1)^{p+q+1}T(\bar{\partial}\varphi )
\,\,\,(T\in {\cal D}^{p,q}(M),\varphi \in A^{n-p,n-q}(M))
\] 
We define the exterior derivative $d$ by 
\[
d := \partial + \bar{\partial}.
\]
\end{definition} 
\begin{definition} 
$T \in {\cal D}^{p,q}(M)$ is said to be closed, if 
$dT = 0$ holds. 
A $(p,p)$-current $T$ is real in case $T = \bar{T}$ in the sense
that $\overline{T(\varphi)} = T(\bar{\varphi})$ holds 
for all $\varphi\in A_{c}^{n-p,n-p}(M)$. 
A real $(p,p)$-current $T$ on $M$ is said to be positive, if
\[
(\sqrt{-1})^{\frac{p(p-1)}{2}}T(\eta\wedge \bar{\eta}) \geq 0
\]
holds for every  $\eta\in A^{n-p}_{c}(M)$. 
\end{definition}
The above definition of positivity of currents is somewhat misleading 
for algebraic geometers. 
It might be appropriate to say {\bf  pseudoeffective currents} 
instead of positive currents. 
\begin{example}
Let $V$ be a subvariety of codimension $p$ in $M$. 
Then $V$ is a closed positive $(p,p)$-current 
on $M$ by 
\[
V(\varphi ):= \int_{V_{reg}}\varphi \,\,\, (\mbox{for}\,\,\varphi\in A^{n-p,n-p}_{c}(M)).
\]
\end{example}
\begin{example} 
Let $\phi$ be a $C^{\infty}$-closed $(p,q)$-form on $M$. 
Then $\phi$ is a closed $(p,q)$-current on $M$ by 
\[
\phi (\varphi ) := \int_{M}\phi\wedge \varphi \,\,\,
(\mbox{for}\,\,\varphi \in A_{c}^{n-p,n-q}(M)) 
\]
\end{example}
\subsection{Multiplier ideal sheaves}
In this subsection $L$ will denote a holomorphic line bundle on a complex manifold $M$. 
\begin{definition}
A  singular hermitian metric $h$ on $L$ is given by
\[
h = e^{-\varphi}\cdot h_{0},
\]
where $h_{0}$ is a $C^{\infty}$-hermitian metric on $L$ and 
$\varphi\in L^{1}_{loc}(M)$ is an arbitrary function on $M$.
We call $\varphi$ a  weight function of $h$.
\end{definition}
The curvature current $\Theta_{h}$ of the singular hermitian line
bundle $(L,h)$ is defined by
\[
\Theta_{h} := \Theta_{h_{0}} + \sqrt{-1}\partial\bar{\partial}\varphi ,
\]
where $\partial\bar{\partial}$ is taken in the sense of a current.
The $L^{2}$-sheaf ${\cal L}^{2}(L,h)$ of the singular hermitian
line bundle $(L,h)$ is defined by
\[
{\cal L}^{2}(L,h) := \{ \sigma\in\Gamma (U,{\cal O}_{M}(L))\mid 
\, h(\sigma ,\sigma )\in L^{1}_{loc}(U)\} ,
\]
where $U$ runs over the  open subsets of $M$.
In this case there exists an ideal sheaf ${\cal I}(h)$ such that
\[
{\cal L}^{2}(L,h) = {\cal O}_{M}(L)\otimes {\cal I}(h)
\]
holds.  We call ${\cal I}(h)$ the {\bf multiplier ideal sheaf} of $(L,h)$.
If we write $h$ as 
\[
h = e^{-\varphi}\cdot h_{0},
\]
where $h_{0}$ is a $C^{\infty}$ hermitian metric on $L$ and 
$\varphi\in L^{1}_{loc}(M)$ is the weight function, we see that
\[
{\cal I}(h) = {\cal L}^{2}({\cal O}_{M},e^{-\varphi})
\]
holds.
Also we define 
\[
{\cal I}_{\infty}(h) = {\cal L}^{\infty}({\cal O}_{M},e^{-\varphi})
\]
and call it the $L^{\infty}$-{\bf multiplier ideal sheaf} of $(L,h)$.

Let $D$ be an effective {\bf  R}-divisor on $M$ and let 
\[
\sum_{i} a_{i}D_{i}
\]
be the irreducible decomposition of $D$.
Let $\sigma_{i}$ be a global section of ${\cal O}_{M}(D_{i})$
with divisor $D_{i}$. 
Let $h_{i}$ be a $C^{\infty}$-hermitian metric on ${\cal O}_{M}(D_{i})$.
Then 
\[
h = \frac{\prod_{i}h_{i}^{a_{i}}}{\prod_{i}h_{i}(\sigma_{i},\sigma_{i})^{a_{i}}}
\]
is a singular hermitian metric on the {\bf  R}-line bundle
${\cal O}_{M}(D)$.
It is clear that $h$ is independent of the choice of 
$h_{i}$'s.  
We define the multiplier sheaf ${\cal I}(D)$ associated with $D$ by
\[
{\cal I}(D) := {\cal I}(h) 
= {\cal L}^{2}({\cal O}_{X},\frac{1}{\prod_{i}h_{i}(\sigma_{i},\sigma_{i})^{a_{i}}}).
\]
If $\mbox{Supp}\,D$ is a divisor with normal crossings,
\[
{\cal I}(D) = {\cal O}_{M}(-[D])
\]
holds, where $[D] := \sum_{i}[a_{i}]D_{i}$ (for a real 
number $a$, $[a]$ 
denotes the largest integer smaller than or equal to $a$).

The following terminology is fundamental in this paper. 
\begin{definition}
$L$ is said to be pseudoeffective, if there exists 
a singular hermitian metric $h$ on $L$ such that 
the curvature current 
$\Theta_{h}$ is a closed positive current.

Also a singular hermitian line bundle $(L,h)$ is said to be pseudoeffective, 
if the curvature current $\Theta_{h}$ is a closed positive current.
\end{definition}
It is easy to see that a line bundle $L$ on a smooth projective 
manifold $M$ is 
pseudoeffective,
if and only if for an ample line bundle $H$ on $M$,  
$L + \epsilon H$ is {\bf Q}-effective (or big) 
for every positive rational number $\epsilon$ (cf. \cite{d}).

If $\{\sigma_{i}\}$ are a finite number of global holomorphic sections of $L$,
for every positive rational number $\alpha$ and a $C^{\infty}$-function 
$\phi$,
\[
h := e^{-\phi}\cdot\frac{h_{0}^{\alpha}}{(\sum_{i}h_{0}(\sigma_{i},\sigma_{i}))^{\alpha}}
\]
defines a singular hermitian metric  on 
$\alpha L$, where $h_{0}$ is a $C^{\infty}$-hermitian metric on $L$
(note that the righthandside is independent of $h_{0}$).
We call such a metric $h$ a singular hermitian metric 
on $\alpha L$ with  {\bf algebraic singularities}.
Singular hermitian metrics with algebraic singularities 
are particulary easy to handle, because its multiplier 
ideal sheaf or that of the multiple of the metric can 
be controlled by taking  suitable successive blowing ups 
such that the total transform of the divisor
$\sum_{i}(\sigma_{i})$ is a divisor with normal crossings. 

By  definition  a multiplier ideal sheaf has the following
property which will be used later.
\begin{lemma}
Let $(L,h)$ be a singular hermitian line bundle on a complex 
manifold $M$ such that $\Theta_{h}$ is bounded from below by 
a $C^{\infty}$-$(1,1)$-form.
Let $f : N\longrightarrow M$ be a modification. 
Then $(f^{*}L,f^{*}h)$ is a singular hermitian line bundle on $N$
and 
\[
f_{*}{\cal I}(f^{*}h) \subseteq {\cal I}(h)
\]
holds.
\end{lemma}
{\bf Proof.}
First we note that $f_{*}(f^{*}L) = L$ holds. 
Let $x\in M$ be an arbitrary point of $M$. 
Let $U$ be a neighbourhood of $x$ and let 
$\sigma$ be a holomorphic section of $L$ on $U$ such that 
\[
\int_{f^{-1}(U)}f^{*}h(\sigma ,\sigma )\, dV_{N} < \infty
\]
holds, where $dV_{N}$ denote a $C^{\infty}$ volume form 
on $N$.
Let $dV_{M}$ be a $C^{\infty}$-volume form on $M$.
Then if we shrink $U$ a little bit, we may assume that 
there exists a positive constant $C$ such that 
\[
f^{*}dV_{M} \leq C\cdot dV_{N}
\]
holds on $f^{-1}(U)$. 
Hence we see that 
\[
\int_{U}h(\sigma ,\sigma )\, dV_{M} < \infty
\]
holds. 
\vspace{5mm} {\bf Q.E.D.} \\ 
The following theorem is fundamental in the applications 
of multiplier ideal sheaves. 
\begin{theorem}(Nadel's vanishing theorem \cite[p.561]{n})
Let $(L,h)$ be a singular hermitian line bundle on a compact K\"{a}hler
manifold $M$ and let $\omega$ be a K\"{a}hler form on $M$.
Suppose that $\Theta_{h}$ is strictly positive, i.e., there exists
a positive constant $\varepsilon$ such that
\[
\Theta_{h} \geq \varepsilon\omega
\]
holds.
Then ${\cal I}(h)$ is a coherent sheaf of ${\cal O}_{M}$ ideal 
and for every $q\geq 1$
\[
H^{q}(M,{\cal O}_{M}(K_{M}+L)\otimes{\cal I}(h)) = 0
\]
holds.
\end{theorem}

We note that the multiplier ideal sheaf of a singular hermitian {\bf R}-line 
bundle is well defined because the multiplier ideal sheaf is defined 
in terms of the weight function.
Sometimes it is useful to consider the 
following variant of multiplier ideal sheaves. 
\begin{definition}
Let $h_{L}$ be a singular hermitian metric on a line bundle $L$.
Suppose that the curvature of $h_{L}$ is a positive current on $X$.
We set 
\[
\bar{\cal I}(h_{L}) 
:= \lim_{\varepsilon\downarrow 0}{\cal I}(h_{L}^{1+\varepsilon})
\]
and call it   the closure of ${\cal I}(h_{L})$. 
\end{definition} 
As you see later, the closure of a multiplier ideal  sheaf is easier to handle 
than the original multiplier ideal sheaf in some respect. 

Next we shall consider the restriction of singular hermitian 
line bundles to subvarieties. 
\begin{definition}
Let $h$ be a singular hermitian metric on $L$ given by 
\[
h = e^{-\varphi}\cdot h_{0},
\]
where $h_{0}$ is a $C^{\infty}$-hermitian metric on $L$ and 
$\varphi\in L^{1}_{loc}(M)$ is an uppersemicontinuous function.
Here $L^{1}_{loc}(M)$ denotes the set of locally integrable functions 
(not the set of classes of almost everywhere equal locally 
integrable functions on $M$). 

For a subvariety $V$ of $M$, we say that the restriction 
$h\mid_{V}$ is well defined, if 
$\varphi$ is not identically $-\infty$ on $V$. 
\end{definition}
Let $(L,h)$,$h_{0}$,$V$, $\varphi$ be as in Definition 2.6.
Suppose that the curvature current $\Theta_{h}$ is bounded 
from below by some $C^{\infty}$-(1,1)-form. 
Then $\varphi$ is an almost plurisubharmonic 
function, i.e. locally a sum of a plurisubharmonic function and a $C^{\infty}$-function.
Let $\pi :\tilde{V} \longrightarrow V$ be an arbitrary 
resolution of $V$. 
Then  $\pi^{*}(\varphi\mid_{V})$ is locally integrable on $\tilde{V}$, 
since $\varphi$ is almost plurisubharmonic. 
Hence 
\[
\pi^{*}(\Theta_{h}\mid_{V}) := \Theta_{\pi^{*}h_{0}\mid_{V}}
+ \sqrt{-1}\partial\bar{\partial}\pi^{*}(\varphi\mid_{V})
\]
is well defined.

\begin{definition}
Let $\varphi$ be a plurisubharmonic function on a unit 
open polydisk  $\Delta^{n}$ with center $O$.
We define the Lelong number of $\varphi$ at $O$ by 
\[
\nu (\varphi ,O) := \liminf_{x\rightarrow O}\frac{\varphi (x)}{\log \mid x\mid},
\]
where $\mid x\mid  = (\sum\mid x_{i}\mid^{2})^{1/2}$.
Let $T$ be a closed positive $(1,1)$-current on 
a unit open polydisk $\Delta^{n}$.
Then by $\partial\bar{\partial}$-Poincar\'{e} lemma
there exists a plurisubharmonic function  $\phi$ 
on $\Delta^{n}$ such that
\[
T = \frac{\sqrt{-1}}{\pi}\partial\bar{\partial}\phi .
\]
We define the Lelong number $\nu (T,O)$ at $O$ by
\[
\nu (T,O) := \nu (\phi ,O).
\]
It is easy to see that $\nu (T,O)$ is independent of the choice of
$\phi$ and local coordinates around $O$.
For an analytic subset $V$ of a complex manifold $X$, we set 
\[
\nu (T,V) = \inf_{x\in V}\nu (T,x).
\]
\end{definition}
\begin{remark} More generally 
the Lelong number is defined for a closed positive
$(k,k)$-current on a complex manifold.
\end{remark}

\begin{theorem}(\cite[p.53, Main Theorem]{s})
Let $T$ be a closed positive $(k,k)$-current on a complex manifold
$M$.
Then for every $c > 0$
\[
\{ x\in M\mid \nu (T,x)\geq c\}
\]
is a subvariety of codimension $\geq k$
in $M$.
\end{theorem}

The following lemma shows a rough relationship between 
the Lelong number of $\nu(\Theta_{h},x)$ at $x\in X$ and the stalk of the multiplier
ideal sheaf ${\cal I}(h)_{x}$ at $x$. 

\begin{lemma}(\cite[p.284, Lemma 7]{b}\cite{b2},\cite[p.85, Lemma 5.3]{s})
Let $\varphi$ be a plurisubharmonic function on 
the open unit polydisk $\Delta^{n}$ with center $O$.
Suppose that $e^{-\varphi}$ is not locally integrable 
around $O$.
Then we have that
\[
\nu (\varphi ,O)\geq 2
\]
holds.
And if
\[
\nu (\varphi ,O) > 2n
\]
holds, then $e^{-\varphi}$ is not locally integrable around $O$.
\end{lemma}
Let $(L,h)$ be a pseudoeffective singular hermitian line bundle 
on a complex manifold $M$. 
The {\bf closure} $\bar{\cal I}(h)$ of the multiplier ideal 
sheaf ${\cal I}(h)$ can be analysed in terms of Lelong numbers
in the following way.  
We note that $\bar{\cal I}(h)$ is coherent ideal sheaf on $M$ by Theorem 2.1.

In the case of $\dim M = 1$, we can compute $\bar{\cal I}(h)$ in terms of the Lelong number $\nu (\Theta_{h},x) (x\in M)$.  
In fact in this case $\bar{\cal I}(h)$ is locally free and 
\[
\bar{\cal I}(h) = {\cal O}_{M}(-\sum_{x\in M}[\nu (\Theta_{h},x)]x)
\]
holds by Lemma 2.2, because $2  = 2\dim M$. 

In the case of $\dim M \geq 2$,
let $f : N \longrightarrow M$ be a modification such that 
$f^{*}\bar{\cal I}(h)$ is locally free. 
If we take $f$ properly, we may assume that there exists a divisor $F = \sum_{i} F_{i}$ with normal crossings on $Y$ such that 
\[
K_{N} = f^{*}K_{M} + \sum_{i} a_{i}F_{i}
\]
and 
\[
\bar{\cal I}(h) = f_{*}{\cal O}_{N}(-\sum_{i} b_{i}F_{i})
\]
hold on $Y$ for some nonnegative integers $\{ a_{i}\}$ and $\{ b_{i}\}$. 
Let $y \in F_{i}- \sum_{j\neq i}F_{j}$ and let 
$(U,z_{1},\ldots ,z_{n})$ be a local corrdinate neighbourhood  of $y$ which is biholomorphic to the open unit disk $\Delta^{n}$ with center $O$ in $\mbox{\bf C}^{n}(n = \dim M)$ and 
\[
U \cap F_{i} = \{ p\in U\mid z_{1}(p) = 0\}
\]
holds.
For $q\in \Delta^{n-1}$, we  set $\Delta (q):= \{ p\in U\mid 
(z_{2}(p),\ldots ,z_{n}(p)) = q\}$.
Then considering the family of the restriction $\{ \Theta_{h}\mid_{\Delta (q)}\}$ for very general $q\in \Delta^{n-1}$, by Lemma 2.2, we see that   
\[
b_{i} = [\nu (f^{*}\Theta_{h},F_{i}) - a_{i}]
\]
holds for every $i$. 
In this way $\bar{\cal I}(h)$ is determined by the {\bf Lelong numbers} of 
the curvature current on some modification. 
This is not the case, unless we take the closure as in the following example. 
\begin{example}
Let $h_{P}$ be a singular hermitian metric on the trivial line bundle on 
the open unit polydisk $\Delta$ with center $O$ in {\bf C} 
defined by
\[
h_{P} = \frac{\parallel\,\cdot\,\parallel^{2}}{\mid z\mid^{2}(\log \mid z\mid  )^{2}}.
\]
Then $\nu (\Theta_{h_{P}},0) = 1$ holds. 
But ${\cal I}(h_{P}) = {\cal O}_{\Delta}$ holds. 
On the other hand $\bar{\cal I}(h_{P}) = {\cal M}_{0}$ holds, 
where ${\cal M}_{0}$ is the ideal sheaf of $0\in \Delta$.  
\end{example}

\subsection{Analytic Zariski decompositions}
In this subsection we shall introduce the notion of analytic Zariski decompositions. 
By using analytic Zariski decompositions, we can handle  big line bundles
like  nef and big line bundles.
\begin{definition}
Let $M$ be a compact complex manifold and let $L$ be a holomorphic line bundle
on $M$.  A singular hermitian metric $h$ on $L$ is said to be 
an analytic Zariski decomposition, if the followings hold.
\begin{enumerate}
\item $\Theta_{h}$ is a closed positive current,
\item for every $m\geq 0$, the natural inclusion
\[
H^{0}(M,{\cal O}_{M}(mL)\otimes{\cal I}(h^{m}))\rightarrow
H^{0}(M,{\cal O}_{M}(mL))
\]
is an isomorphim.
\end{enumerate}
\end{definition}
\begin{remark} If an AZD exists on a line bundle $L$ on a smooth projective
variety $M$, $L$ is pseudoeffective by the condition 1 above.
\end{remark}

\begin{theorem}(\cite{tu,tu2})
 Let $L$ be a big line  bundle on a smooth projective variety
$M$.  Then $L$ has an AZD. 
\end{theorem}
As for the existence for general pseudoeffective line bundles, 
now we have the following theorem.
\begin{theorem}(\cite{d-p-s})
Let $X$ be a smooth projective variety and let $L$ be a pseudoeffective 
line bundle on $X$.  Then $L$ has an AZD.
\end{theorem}
{\bf Proof of Theorem 2.4}. Let  $h_{0}$ be a fixed $C^{\infty}$-hermitian metric on $L$.
Let $E$ be the set of singular hermitian metric on $L$ defined by
\[
E = \{ h ; h : \mbox{lowersemicontinuous singular hermitian metric on $L$}, 
\]
\[
\hspace{70mm}\Theta_{h}\,
\mbox{is positive}, \frac{h}{h_{0}}\geq 1 \}.
\]
Since $L$ is pseudoeffective, $E$ is nonempty.
We set 
\[
h_{L} = h_{0}\cdot\inf_{h\in E}\frac{h}{h_{0}},
\]
where the infimum is taken pointwise. 
The supremum of a family of plurisubharmonic functions 
uniformly bounded from above is known to be again plurisubharmonic, 
if we modify the supremum on a set of measure $0$(i.e., if we take the uppersemicontinuous envelope) by the following theorem of P. Lelong.

\begin{theorem}(\cite[p.26, Theorem 5]{l})
Let $\{\varphi_{t}\}_{t\in T}$ be a family of plurisubharmonic functions  
on a domain $\Omega$ 
which is uniformly bounded from above on every compact subset of $\Omega$.
Then $\psi = \sup_{t\in T}\varphi_{t}$ has a minimum 
uppersemicontinuous majorant $\psi^{*}$ which is plurisubharmonic.
\end{theorem}
\begin{remark} In the above theorem the equality 
$\psi = \psi^{*}$ holds outside of a set of measure $0$(cf.\cite[p.29]{l}). 
\end{remark}

By Theorem 2.5 we see that $h_{L}$ is also a 
singular hermitian metric on $L$ with $\Theta_{h}\geq 0$.
Suppose that there exists a nontrivial section 
$\sigma\in \Gamma (X,{\cal O}_{X}(mL))$ for some $m$ (otherwise the 
second condition in Definition 3.1 is empty).
We note that  
\[
\frac{1}{\mid\sigma\mid^{\frac{2}{m}}} 
\]
gives the weihgt of a singular hermitian metric on $L$ with curvature 
$2\pi m^{-1}(\sigma )$, where $(\sigma )$ is the current of integration
along the zero set of $\sigma$. 
By the construction we see that there exists a positive constant 
$c$ such that  
\[
\frac{h_{0}}{\mid\sigma\mid^{\frac{2}{m}}} \geq c\cdot h_{L}
\]
holds. 
Hence
\[
\sigma \in H^{0}(X,{\cal O}_{X}(mL)\otimes{\cal I}(h_{L}^{m}))
\]
holds.   This means that $h_{L}$ is an AZD of $L$. 
\vspace{10mm} {\bf Q.E.D.}  \\
The following proposition implies that the multiplier ideal sheaves 
of $h_{L}^{m}(m\geq 1)$ constructed in the proof of
 Theorem 2.4 are independent of 
the choice of the $C^{\infty}$-hermitian metric $h_{0}$.
The proof is trivial.  Hence we omit it. 
\begin{proposition}
$h_{0},h_{0}^{\prime}$ be two $C^{\infty}$-hermitian metrics 
on a pseudoeffective line bundle $L$ on a smooth projective 
variety $X$. 
Let $h_{L},h^{\prime}_{L}$ be the AZD's constructed as in the 
proof of Theorem 2.4 associated with $h_{0},h_{0}^{\prime}$ 
respectively. 
Then 
\[
(\min_{x\in X}\frac{h_{0}}{h_{0}^{\prime}}(x))\cdot h_{L}^{\prime}
     \leq    h_{L} \leq
 (\max_{x\in X}\frac{h_{0}}{h_{0}^{\prime}}(x))\cdot h_{L}^{\prime}
\]
hold.
In particular 
\[
{\cal I}(h_{L}^{m}) = {\cal I}((h_{L}^{\prime})^{m})
\]
holds for every $m\geq 1$.
\end{proposition}
We call the AZD constructed as in the proof of Theorem 2.4  {\bf a canonical 
AZD} of $L$. 
Proposition 2.1 implies that the multiplier ideal sheaves associated with 
the multiples of the canonical AZD are independent of the choice of 
the canonical AZD. 

\subsection{Intersection numbers}
In this subsection we shall define the intersection number 
for a singular hermitian line bundle with positive curvature current 
and an irreducible curve such that the restriction of 
the singular hermitian metric is well defined. 
\begin{definition}
Let $(L,h)$ be a pseudoeffective singular hermitian line bundle on a smooth 
projective variety $X$.
Let $C$ be an irreducible curve on $X$ such that the natural morphism 
${\cal I}(h^{m})\otimes{\cal O}_{C}\rightarrow {\cal O}_{C}$ 
is an isomorphism at the generic point of $C$ for every $m\geq 0$.

The intersection number $(L,h)\cdot C$ is defined by
\[
(L,h)\cdot C := 
\overline{\lim}_{m\rightarrow\infty}m^{-1}\dim H^{0}(C,{\cal O}_{C}(mL)\otimes
{\cal I}(h^{m})/tor),
\]
where $tor$ denotes the torsion part of 
${\cal O}_{C}(mL)\otimes {\cal I}(h^{m})$. 
\end{definition}
If the natural morphism ${\cal I}(h^{m})\otimes{\cal O}_{C}\rightarrow {\cal O}_{C}$ is $0$ at the generic point of $C$ for some $m\geq 1$, to define $(L,h)\cdot C$, 
 $H^{0}(C,{\cal O}_{C}(mL)\otimes
{\cal I}(h^{m})/tor)$ cannot be considered as a subspace of $H^{0}(C,{\cal O}_{C}(mL))$. 
A special important case will be treated in Section 2.5.  

\begin{remark} 
Let $(L,h)$, $C$ be as above. 
Let $\pi : \tilde{C}\longrightarrow C$ be the 
normalization of $C$. 
Then we see that 
\[
(L,h)\cdot C = 
\overline{\lim}_{m\rightarrow\infty}m^{-1}\dim H^{0}(\tilde{C},{\cal O}_{\tilde{C}}
(m\pi^{*}L)\otimes\pi^{*}{\cal I}(h^{m})/tor)
\]
holds. 
This is verified as follows. 
First it is clear that 
\[
(L,h)\cdot C \leq  
\overline{\lim}_{m\rightarrow\infty}m^{-1}\dim H^{0}(\tilde{C},{\cal O}_{\tilde{C}}
(m\pi^{*}L)\otimes\pi^{*}{\cal I}(h^{m})/tor)
\]
holds. 
On the other hand, there exists a nonzero ideal sheaf ${\cal J}$ independent of 
$m\geq 0$ on 
$\tilde{C}$ such that 
\[
H^{0}(\tilde{C},({\cal O}_{\tilde{C}}(m\pi^{*}L)
\otimes\pi^{*}{\cal I}(h^{m})/tor)\otimes{\cal J})
\subseteq 
\pi^{*}H^{0}(C,{\cal O}_{C}(mL)\otimes{\cal I}(h^{m})/tor)
\]
holds. 
For example, we can take ${\cal J}$ to be 
$(\pi^{*}{\cal I}_{\mbox{Sing}(C)})^{r}$ 
where ${\cal I}_{\mbox{Sing}(C)}$ denotes the 
ideal sheaf of the singular locus of $C$ and 
$r$ is a sufficiently large positive integer.
Because $V({\cal I})$ consists of 
a finite number of points, this implies that 
\[
(L,h)\cdot C \geq 
\overline{\lim}_{m\rightarrow\infty}m^{-1}\dim H^{0}(\tilde{C},{\cal O}_{\tilde{C}}
(m\pi^{*}L)\otimes\pi^{*}{\cal I}(h^{m})/tor)
\]
holds. 
The above two inequalities imply the assertion.
\end{remark}
\begin{remark}
Let $(L,h)$, $C$ be as in Definition 2.9. 
We see that 
\[
(L,h)\cdot C = \overline{\lim}_{m\rightarrow\infty}m^{-1}
\dim H^{0}(C,{\cal O}_{C}(mL)\otimes
\bar{{\cal I}}(h^{m})/tor)
\]
always holds.

This can be verified as follows.
First we shall assume that $C$ is smooth. 
By the assumption ${\cal I}(h^{m})/tor$ is an ideal sheaf on $C$.
If 
\[
\deg_{C}{\cal O}_{C}(mL)\otimes{\cal I}(h^{m})/tor > 2g(C)-2
\]
holds, where $g(C)$ denotes the genus of $C$, then 
\[
H^{1}(C,{\cal O}_{C}(mL)\otimes{\cal I}(h^{m})/tor) = 0
\]
holds. 
On the other hand if 
\[
\deg_{C}{\cal O}_{C}(mL)\otimes{\cal I}(h^{m})/tor \leq  2g(C)-2
\]
holds, then there exists a constant $K$ independent of such $m$ such that 
\[
\dim H^{0}(C,{\cal O}_{C}(mL)\otimes{\cal I}(h^{m})/tor) \leq K
\]
holds. 

Hence we see that 
\[
(\flat ) \hspace{10mm} (L,h)\cdot C =  \overline{\lim}_{m\rightarrow\infty}m^{-1}\deg_{C}{\cal O}_{C}(mL)\otimes{\cal I}(h^{m})/tor
\]
holds by the Riemann-Roch theorem. 
By the same reason, we see that 
\[
\overline{\lim}_{m\rightarrow\infty}m^{-1}\dim H^{0}(C,{\cal O}_{C}(mL)\otimes
\bar{{\cal I}}(h^{m})/tor)
= 
\overline{\lim}_{m\rightarrow\infty}m^{-1}\deg_{C}{\cal O}_{C}(mL)\otimes
\bar{{\cal I}}(h^{m})/tor 
\]
holds. 
On the other hand, for every $\epsilon > 0$
\[
\overline{\lim}_{m\rightarrow\infty}m^{-1}\deg_{C}{\cal O}_{C}(mL)\otimes
{\cal I}(h^{m})/tor
\]
\[
\hspace{20mm}
\leq \overline{\lim}_{m\rightarrow\infty}m^{-1} 
\deg_{C}{\cal O}_{C}(\lceil (1+2\epsilon )m\rceil L)\otimes
\bar{{\cal I}}(h^{(1+\epsilon )m})/tor 
\]
holds, since ${\cal I}(h^{(1+2\epsilon )m})
\subseteq \bar{{\cal I}}(h^{(1+\epsilon )m})$
holds for every $m \geq 0$. 
And also 
\[
\lim_{\epsilon\downarrow 0}(\overline{\lim}_{m\rightarrow\infty}m^{-1} 
\deg_{C}{\cal O}_{C}(\lceil (1+2\epsilon )m\rceil L)\otimes
\bar{{\cal I}}(h^{(1+\epsilon )m})/tor)
\]
\[
\hspace{20mm} = 
\lim_{\epsilon\downarrow 0}((\overline{\lim}_{m\rightarrow\infty}m^{-1} 
\deg_{C}{\cal O}_{C}(\lceil (1+\epsilon )m\rceil L)\otimes
\bar{{\cal I}}(h^{(1+\epsilon )m})/tor) + \epsilon\, L\cdot C) 
\] 
\[
\hspace{20mm} = 
\overline{\lim}_{m\rightarrow\infty}m^{-1} 
\deg_{C}{\cal O}_{C}(m L)\otimes
\bar{{\cal I}}(h^{m})/tor 
\] 
hold. 
We note that $\bar{{\cal I}}(h^{m})\subseteq {\cal I}(h^{m})$ holds for 
every $m\geq 0$ by their definitions.
Hence we have that 
\[
\overline{\lim}_{m\rightarrow\infty}m^{-1}\deg_{C} {\cal O}_{C}(mL)\otimes
{\cal I}(h^{m})/tor
= 
\overline{\lim}_{m\rightarrow\infty}m^{-1} 
\deg_{C} {\cal O}_{C}(mL)\otimes
\bar{{\cal I}}(h^{m})/tor 
\]
holds.
By the above argument we see that 
\[
 (L,h)\cdot C = \overline{\lim}_{m\rightarrow\infty}m^{-1}\dim H^{0}(C,{\cal O}_{C}(mL)\otimes
\bar{{\cal I}}(h^{m})/tor)
\]
 holds.

If $C$ is singular, by the argument as in Remark 2.4, 
we can easily deduce the same conclusion by considering 
the normalization $\pi :\tilde{C}\longrightarrow C$.

Since the closure of multiplier a multiplier ideal sheaf is easier to handle 
as you see in this paper, it might be better to use the above formula as 
the definition of the intersection number.  
\end{remark}
Let $(L,h)$ and $C$ be as above.
Assume that $h\mid_{C}$ is well defined. 
Let 
\[
\pi : \tilde{C}\longrightarrow C
\]
be the normalization of $C$. 
We define the multiplier ideal sheaf \\
${\cal I}(h^{m}\mid_{C}) (m\geq 0)$ 
on $C$ by 
\[
{\cal I}(h^{m}\mid_{C}) := \pi_{*}{\cal I}(\pi^{*}h^{m}\mid_{C}). 
\]
We note that ${\cal I}(h^{m}\mid_{C})$ is not necessary 
a subsheaf of ${\cal O}_{C}$, if $C$ is nonnormal. 
And the Lelong number $\nu (\Theta_{h}\mid_{C},x) (x\in C)$ by 
\[
\nu (\Theta_{h},x) = \sum_{\tilde{x}\in \pi^{-1}(x)}
\nu (\pi^{*}\Theta_{h}\mid_{C},\tilde{x}). 
\]
\begin{proposition}
Let $(L,h)$ be a pseudoeffective singular hermitian line bundle on 
a smooth projective variety $X$.
Let $C$ be an irreducible curve on $X$ such that $h\mid_{C}$ 
is well defined. 
Suppose that $(L,h)\cdot C = 0$ holds. 
Then 
\[
\Theta_{h}\mid_{C} = 2\pi\sum_{x\in C}\nu (\Theta_{h}\mid_{C},x)x
\]
holds in the sense that 
\[
\pi^{*}(\Theta_{h}\mid_{C}) = 
2\pi\sum_{\tilde{x}\in \tilde{C}}\nu (\pi^{*}\Theta_{h}\mid_{C},\tilde{x})\tilde{x}
\]
holds. 
\end{proposition}
{\bf Proof of Proposition 2.2.}  
First we quote the following $L^{2}$-extension theorem.  
\begin{theorem}(\cite[p.197, Theorem]{o-t})
Let $\Omega$ be a bounded pseudoconvex domain in ${\bf C}^{n}$,
$\psi : \Omega\longrightarrow {\bf R}\cup \{-\infty\}$ a plurisubharmonic
function and $H\subset {\bf C}^{n}$ a complex hyperplane. 

Then there exists a constant $C$ depending only on the diameter of $\Omega$ 
such that for any holomorphic function $f$ on $\Omega\cap H$ satisfying
\[
\int_{\Omega\cap H}e^{-\psi}\mid f\mid^{2}dV_{n-1} < \infty,
\]
where $dV_{n-1}$ denotes the $(2n-2)$-dimensional Lebesgue measure, 
there exists a holomorphic function $F$ on $\Omega$ satisfying 
$F\mid_{\Omega\cap H} = f$ and 
\[
\int_{\Omega}e^{-\psi}\mid F\mid^{2}dV_{n}
\leq 
C\cdot \int_{\Omega\cap H}e^{-\psi}\mid f\mid^{2}dV_{n-1}.
\]
\end{theorem}
\begin{lemma}
Let $S$ be the singular points of $C$ with reduced structure and 
let ${\cal I}_{S}$ denote the ideal of $S$. 
Then there exists a positive integer $a$ such that 
\[
{\cal I}(h^{m}\mid_{C})\otimes {\cal I}_{S}^{a} \subset {\cal I}(h^{m})\otimes{\cal O}_{C}
\]
hold for every $m$. 
\end{lemma}
{\bf Proof of Lemma 2.3}.
In fact let 
\[
f : \tilde{X}\longrightarrow X
\]
be an embedded resolution of $C$ and let $\tilde{C}$ denote 
the strict transform of $C$ in $\tilde{X}$. 
Since $\tilde{C}$ is locally a smooth complete intersection of 
smooth divisors, for $\tilde{x}\in \tilde{C}$, by
the successive use of Theorem 2.6 every element 
 of ${\cal I}(f^{*}h^{m}\mid_{\tilde{C}})_{\tilde{x}}$
can be extended to an element of ${\cal I}(f^{*}h^{m})_{\tilde{x}}$. 
This means  that 
\[
{\cal I}(f^{*}h^{m}\mid_{\tilde{C}}) \subseteq 
{\cal I}(f^{*}h^{m})\mid_{\tilde{C}}
\]
holds. 
By the definition of ${\cal I}(h^{m}\mid_{C})$ we see that 
\[
{\cal I}(h^{m}\mid_{C}) = f_{*}({\cal I}(f^{*}h^{m}\mid_{\tilde{C}}))
\]
holds. 
Hence we have that 
\[
\mbox{($+$)} \hspace{10mm}
{\cal I}(h^{m}\mid_{C}) \subseteq f_{*}({\cal I}(f^{*}h^{m})\mid_{\tilde{C}})
\]
holds. 

Let $f^{*}(C)$ be the total transform of $C$. 

First we note that if a germ of ${\cal I}(f^{*}h^{m})\mid_{\tilde{C}}$
is identically $0$ along the scheme theoretic intersection 
$(f^{*}(C) - \tilde{C})\cap \tilde{C}$, 
it extends to a germ of ${\cal I}(f^{*}h^{m})\mid_{f^{*}(C)}$ 
by setting identically $0$ on the branches of $f^{*}(C)$ except $\tilde{C}$. 

Next we note that 
\[
f_{*}({\cal I}(f^{*}h^{m})\otimes
{\cal O}_{f^{*}(C)})
\subseteq {\cal I}(h^{m})\otimes{\cal O}_{C}
\]
holds by Lemma 2.1.

By these facts and ($+$), we see that there exists a positive integer $a$ 
independent of $m$ such that 
\[
{\cal I}(h^{m}\mid_{C})\otimes {\cal I}_{S}^{a} \subseteq 
{\cal I}(h^{m})\otimes {\cal O}_{C}
\]
holds. 
This completes the proof of Lemma 2.3. {\bf Q.E.D.}
\vspace{10mm}\\

By \cite[p.111, Lemma 9.5]{s} we see that 
\[
\Theta_{h}\mid_{C}- 2\pi\sum_{x\in C}\nu (\Theta_{h}\mid_{C},x)x
\]
is a positive current on $C$.
Hence 
\[
L\cdot C - \sum_{x\in C}\nu (\Theta_{h}\mid_{C},x)
\]
is a nonnegative number. 
Let $\pi : \tilde{C}\longrightarrow C$ be the 
normalization of $C$. 
Then by Lemma 2.2 and the definition of 
$\nu (\Theta_{h}\mid_{C})$, we see that 
\[
\deg_{\tilde{C}}({\cal O}_{\tilde{C}}(m\pi^{*}L)\otimes{\cal I}(\pi^{*}(h^{m}\mid_{C}))
\geq  (L\cdot C - \sum_{x\in C}\nu (\Theta_{h}\mid_{C},x))m
\]
holds. 
Hence we see that 
\[
\lim_{m\rightarrow\infty}m^{-1}\deg_{\tilde{C}}({\cal O}_{\tilde{C}}(m\pi^{*}L)\otimes{\cal I}(\pi^{*}(h^{m}\mid_{C}))
\geq 
L\cdot C - \sum_{x\in C}\nu (\Theta_{h}\mid_{C},x)\geq 0
\]
hold.
By the Riemann-Roch theorem for curves and the Kodaira vanishing theorem, 
we see that 
if 
\[
L\cdot C - \sum_{x\in C}\nu (\Theta_{h}\mid_{C},x) > 0
\]
holds, then 
\[
\lim_{m\rightarrow\infty}m^{-1}\dim H^{0}(\tilde{C},{\cal O}_{\tilde{C}}(m\pi^{*}L)\otimes{\cal I}(\pi^{*}(h^{m}\mid_{C}))
\geq 
L\cdot C - \sum_{x\in C}\nu (\Theta_{h}\mid_{C},x)
\]
holds. 

By Lemma 2.3, this means that $(L,h)\cdot C$ is always nonnegative and 
\[
(L,h)\cdot C > 0
\]
holds, when
\[
L\cdot C - \sum_{x\in C}\nu (\Theta_{h}\mid_{C},x) > 0
\]
holds.  
Hence if $(L,h)\cdot C = 0$ holds, then
\[
L\cdot C = \sum_{x\in C}\nu (\Theta_{h}\mid_{C},x)
\]
holds. 
This implies that 
\[
\Theta_{h}\mid_{C}  = 2\pi\sum_{x\in C}\nu (\Theta_{h}\mid_{C},x)x
\]
holds.
This completes the proof of Proposition 2.2. {\bf Q.E.D.} 
\begin{definition} 
Let $(L,h)$ be a pseudoeffective singular hermitian line bundle on 
a smooth projective variety $X$. 
$(L,h)$ is said to be numerically trivial, if for every 
irreducible curve $C$ on $X$ such that $h\mid_{C}$ is 
well defined, 
\[
(L,h)\cdot C = 0
\]
holds. 
\end{definition}
\subsection{Restriction of the intersection theory to 
divisors}
In the previous subsection we define an 
intersection number of a singular hermitian line bundle 
with positive curvature and an irreducible curve on 
which the restriction of the singular hermitian metric 
is well defined.  In this subsection we shall consider 
the case that the restriction of the singular hermitian 
metric is not well defined. 

Let $(L,h)$ be a pseudoeffective singular hermitian line bundle on 
a smooth projective variety $X$. 

Let $D$ be a smooth divisor on $X$. 
We set  
\[
v_{m}(D) = \mbox{mult}_{D}\mbox{Spec}({\cal O}_{X}/{\cal I}(h^{m}))
\]
and 
\[
\tilde{\cal I}_{D}(h^{m}) = {\cal O}_{D}(v_{m}(D)D)\otimes {\cal I}(h^{m}).
\] 
Then $\tilde{\cal I}_{D}(h^{m})$ is an ideal sheaf on $D$ (it is torsion free, 
since $D$ is smooth). 

Let $x\in D$ be an arbitrary point of $D$ and let 
$(U,z_{1},\ldots ,z_{n}) (n:= \dim X)$ be a local coordinate neighbourhood of 
$x$ which is 
biholomorphic to the unit open polydisk $\Delta^{n}$ with center $O$ in 
$\mbox{\bf C}^{n}$ and 
\[
U \cap D = \{ p\in U\mid z_{1}(p) = 0\}
\] 
holds. 
For $q\in \Delta^{n-1}$, we  set $\Delta (q):= \{ p\in U\mid 
(z_{2}(p),\ldots ,z_{n}(p)) = q\}$.
Then considering the family of the restriction $\{ \Theta_{h}\mid_{\Delta (q)}\}$ for very general $q\in \Delta^{n-1}$, by Lemma 2.2, we see that  
\[
m\cdot\nu (\Theta_{h},D) -1 \leq v_{m}(D) 
\leq m\cdot\nu (\Theta_{h},D)
\]
holds. 

We define the ideal sheaves $\sqrt[m]{\tilde{\cal I}_{D}(h^{m})}$ 
on $D$ by 
\[
\sqrt[m]{\tilde{\cal I}_{D}(h^{m})}_{x}
:= \cup {\cal I}(\frac{1}{m}(\sigma))_{x} (x\in D),
\]
where $\sigma$ runs all the germs of $\tilde{\cal I}_{D}(h^{m})_{x}$.
And we set 
\[
{\cal I}_{D}(h):= \cap_{m\geq 1}\sqrt[m]{\tilde{\cal I}_{D}(h^{m})}
\]
and call it {\bf the multipler ideal of $h$ on $D$}.
Also we set 
\[
\bar{\cal I}_{D}(h):= \lim_{\varepsilon\downarrow 0}
{\cal I}_{D}(h^{1+\varepsilon}).
\]
See Theorem 2.8 below for the reason why we define ${\cal I}_{D}(h)$ in this way. 

Let $C$ be an irreducible curve in $D$ such that 
the natural morphism 
\[
\tilde{\cal I}_{D}(h^{m})\otimes {\cal O}_{C}\rightarrow {\cal O}_{C}
\]
is an isomorphism at the generic point of $C$ for every $m\geq 0$. 
In this case we can define the intersection number $(L,h)\cdot C$ by  
\[
(L,h)\cdot C := 
\overline{\lim}_{m\rightarrow\infty} 
m^{-1}\dim H^{0}(C,{\cal O}_{C}(mL-v_{m}(D)D)
\otimes\tilde{\cal I}_{D}(h^{m})/tor).
\]
Then as the formla $(\flat )$ in Remark 2.5, we see that 
\[
(\sharp )\hspace{20mm} (L,h)\cdot C =  (L- \nu (\Theta_{h},D)D)\cdot C + 
\overline{\lim}_{m\rightarrow\infty}
m^{-1}\deg_{C}\tilde{\cal I}_{D}(h^{m})\otimes{\cal O}_{C}
\]
holds. 

We may define the {\bf Lelong number} $\nu_{D}(\Theta_{h},x) (x\in D)$ 
by 
\[
\nu_{D}(\Theta_{h},x) := \overline{\lim}_{m\rightarrow\infty}
m^{-1}\mbox{mult}_{x}\mbox{Spec}({\cal O}_{D}/\tilde{\cal I}_{D}(h^{m})).
\]
Then we see that the set 
\[
S_{D}:= \{ x\in D\mid \nu (\Theta_{h}\mid_{D},x) > 0\}
\]
consists of a countable union of subvarieties on $D$. 
This follows from the approximation theorem \cite[p.380, Proposition 3.7]{d}.
\subsection{Another definition of the intersection numbers}
Let $(L,h)$ be a pseudoeffective singular hermitian line bundle on 
a smooth projective variety $X$. 
And let $C$ be an irreducible curve on $X$ such that 
the restriction $h\mid_{C}$ is well defined. 
Another candidate for the intersection number of $(L,h)$ 
and $C$ is :
\[
(L,h)*C := L\cdot C - \sum_{x\in C}\nu (\Theta_{h}\mid_{C},x).
\]
But we have the following theorem.
\begin{theorem}
With the above notations 
\[
(L,h)\cdot C = (L,h)*C
\]
holds. 
\end{theorem}
{\bf Proof of Theorem 2.7}. 
To prove Theorem 2.7, by taking an embedded resolution of $C$, we may assume that $C$ is smooth, since the intersection number $(L,h)\cdot C$ is defined by using the asymptotics of the  dimension of sections
and $\nu (\Theta_{h}\mid_{C},x) (x\in C)$ is 
defined in terms of the normalization of $C$. 

In fact let $\varpi : Y\longrightarrow X$ be an embedded 
resolution of $C$ and let $\tilde{C}$ be the 
strict transform of $C$. 
Since by the definition of multiplier ideal sheaves
\[
{\cal O}_{X}(K_{X})\otimes {\cal I}(h^{m})
= \varpi_{*}({\cal O}_{Y}(K_{Y})\otimes {\cal I}(\varpi^{*}h^{m}))
\]
holds for every $m\geq 0$, we have that 
\[
\varpi^{*}{\cal I}(h^{m})\otimes 
{\cal O}_{Y}(\varpi^{*}K_{X} - K_{Y})\subseteq 
{\cal I}(\varpi^{*}h^{m}) 
\]
holds for every $m\geq 0$. 
Then by Remark 2.4, we have that 
\[
(L,h)\cdot C 
\leq  \overline{\lim}_{m\rightarrow\infty}
m^{-1}\dim H^{0}(\tilde{C},{\cal O}_{\tilde{C}}
(m\varpi^{*}L)\otimes {\cal I}(\varpi^{*}h^{m})/tor)
\]
holds.
Suppose that 
\[
(\varpi^{*}L,\varpi^{*}h)\cdot \tilde{C} 
= (\varpi^{*}L,\varpi^{*}h)*\tilde{C}
\]
holds. 
Then by the above inequality, we have that 
\[
(L,h)\cdot C \leq (\varpi^{*}L,\varpi^{*}h)*\tilde{C} 
\]
holds. 
Since by definition 
\[
(L,h)*C = (\varpi^{*}L,\varpi^{*}h)*\tilde{C} 
\]
holds, we have that 
\[
(L,h)\cdot C \leq (L,h)*C
\]
holds. 
On the other hand by Lemma 2.3 and Lemma 2.2 (see also the explanataion right
after Lemma 2.2), we see that 
the opposite inequality :
\[
(L,h)\cdot C \geq (L,h)*C
\]
holds. 
Hence we conclude that 
\[
(L,h)\cdot C = (L,h)*C
\]
holds. 

Hereafter we shall assume that $C$ is smooth. 
First we note that for every ample line bundle $H$ on $X$
and a $C^{\infty}$-hermitian metric $h_{H}$ on $H$ with 
strictly positive curvature 
\[
(L\otimes H,h\cdot h_{H})\cdot C = H\cdot C + (L,h)\cdot C
\]
holds by the formula $(\sharp )$ and 
\[
(L\otimes H,h\cdot h_{H})* C = H\cdot C + (L,h)*C
\]
hold.  
Hence we may assume that $h$ is strictly positive. 

Since we already have the inequality :
\[
(L,h)\cdot C \geq (L,h)*C
\]
as above, we only have to show the opposite inequality 
\[
(L,h)\cdot C \leq (L,h)*C
\]
holds. 

First we shall consider the case that 
$h$ has algebraic singularities. 
In this case by taking a suitable modification 
\[
f : \tilde{X}\rightarrow X
\]
we see that there exists an effective {\bf Q}-divisor $D$ 
with normal crossings on $\tilde{X}$ such that 
\[
{\cal I}(f^{*}h^{m}) = {\cal O}_{\tilde{X}}(\lceil -mD\rceil)
\]
holds for every $m\geq 0$, where $\lceil\,\,\,\,\rceil$
denotes the round up.
Let $\tilde{C}$ denote the strict transform of $C$. 
We may assume that $\tilde{C}$ is smooth. 
By this  
\[
\deg_{C}({\cal I}(h^{m})\mid_{C}) 
= - [mD]\cdot\tilde{C} + (K_{\tilde{X}} - f^{*}K_{X})\cdot\tilde{C}
\] 
holds.
On the other hand 
\[
\deg_{C} {\cal I}(h^{m}\mid_{C})
= - \deg_{\tilde{C}} [mD\mid_{\tilde{C}}] 
\]
holds. 
Then since 
\[
\lim_{m\rightarrow\infty}\frac{1}{m}[mD]\cdot\tilde{C} 
= \lim_{m\rightarrow\infty}\frac{1}{m}\deg_{C} [mD\mid_{\tilde{C}}] 
\]
holds, we have that 
\[
\lim_{m\rightarrow \infty}\frac{1}{m}\deg_{C} ({\cal I}(h^{m})\mid_{C}) = 
\lim_{m\rightarrow \infty}\frac{1}{m}\deg_{C} {\cal I}(h^{m}\mid_{C})
\]
holds. 
The lefthandside is equal to 
\[
L\cdot C - (L,h)\cdot C
\]
by the argument in Remark 2.5 (especially by the formula $(\flat )$) and the righthandside is equal to 
\[
L\cdot C - (L,h)*C
\]
by Lemma 2.2. 
Hence if $h$ has  algebraic singularities,
\[
(L,h)*C = (L,h)\cdot C
\]
holds. 

On the other hand by (a slight generalization of) the approximation theorem of \cite[p.380, Proposition 3.7]{d}, 
there exists a sequence of singular hermitian metrics $\{ h_{j}\}_{j=1}^{\infty}$ 
satisfying the following {\bf 6-conditions} :  
\begin{enumerate}
\item $\Theta_{h_{j}}$ is positive for every $j$,
\item 
$\lim_{j\rightarrow\infty} h_{j} = h$ 
holds in the sense of  the convergence of the weight functions as currents on $M$ and $C$,
\item $h_{j}$ has algebraic singularities, 
\item  ${\cal I}(h^{jm}) \subseteq {\cal I}(h_{j}^{jm})$, 
holds for every $m\geq 0$ and $j\geq 1$,
\item $\lim_{j\rightarrow\infty}\bar{\cal I}(h_{j}^{m}) = \bar{\cal I}(h^{m})$
holds for every $m$,
\item $\lim_{j\rightarrow\infty}\bar{\cal I}(h_{j}^{m}\mid_{C}) = \bar{\cal I}(h^{m}\mid_{C})$ holds for every $m$.
\end{enumerate}
The third condition looks a little bit different from 
\cite[Proposition 3.7]{d}.
But it is essentially the same by Lemma 2.2 and the construction of $\{ h_{j}\}$ below. 
The 4-th condition cannot be deduced directly by 
the approximation theorem of \cite[p.380, Proposition 3.7]{d}. 

Let us briefly show how to  construct $\{ h_{j}\}$. 
The following argument is a slight modification of that in \cite{d}. 
First we shall consider the local approximation of a plurisubharmonic function 
by a sequence of plurisubharmonic functions with algebraic singularities. 

Let $\varphi$ be a plurisubharmonic function on $\Delta^{n}$. 
Let $C = \{p \in\Delta^{n}\mid z_{2}(p) = \cdots z_{n}(p) = 0\}$. 
Suppose that $\varphi$ is not identically $-\infty$ on $C$. 
That is to say we are considering the case that $h = e^{-\varphi}$ and $C$ is a smooth curve 
in $\Delta^{n}$. Let $m$ be a positive integer. 
Let ${\cal H}(j\varphi )_{C}$ be the Hilbert space defined by
\[
{\cal H}(j\varphi )_{C}: = \{ f\in {\cal O}(\Delta^{n})\mid 
\int_{\Delta^{n}}\mid f\mid^{2}e^{-j\varphi}d\lambda < \infty \,\,
\mbox{and}\,\,
\int_{C}\mid f\mid^{2}e^{-j\varphi}d\lambda_{C} < \infty \} 
\]
with the inner product 
\[
(f,g) := \frac{1}{2}\int_{\Delta^{n}} f\cdot \bar{g}\cdot e^{-j\varphi}d\lambda
+ \frac{1}{2}\int_{C} f\cdot\bar{g}\cdot e^{-j\varphi}d\lambda_{C}
\]
where $d\lambda$ and $d\lambda_{C}$ is the usual Lebesgue measure on 
$\Delta^{n}$ and $C$ respectively.
Let  $\{ \sigma_{\ell}\}$ be an orthonormal basis of 
 ${\cal H}(j\varphi )_{C}$ and let 
\[
\varphi_{j}:= \frac{1}{2j}\log \sum\mid\sigma_{\ell}\mid^{2}.
\]
Let $\psi$ is the plurisubharmonic function on $\Delta^{n}$ defined by
\[
\psi = (n-1)\log (\sum_{i=2}^{n}\mid z_{i}\mid^{2}).
\]
\begin{proposition}
There exist positive constants $K_{1},K_{2} > 0$ independent of $m$ such that 
\begin{enumerate}
\item  \[
\varphi (z)- \frac{K_{1}}{j} \leq \varphi_{j}(z) \leq 
\sup_{\mid\zeta -z\mid < r}\varphi (\zeta ) +\frac{1}{j}\log (\frac{K_{2}}{r^{n}})
\]
holds for every $z\in C$ and $r < d(z,\partial\Delta^{n})$ 
and 
\[
 \varphi (z)+\frac{1}{2j}\psi (z)- \frac{K_{1}}{j} \leq \varphi_{j}(z) \leq 
\sup_{\mid\zeta -z\mid < r}\varphi (\zeta ) +\frac{1}{j}\log (\frac{K_{2}}{r^{n}})
\]
holds for every $z\in \Delta^{n}- C$ and $r < d(z,\partial\Delta^{n})$,
\item $\nu (\varphi ,z) -n/j\leq \nu (\varphi_{j},z) \leq 
\nu (\varphi ,z)$ holds for every $z\in \Delta^{n}$. 
\end{enumerate}
\end{proposition}
{\bf Proof of Proposition 2.3.} 
We note that 
\[
\varphi_{j}(z) = \sup_{f\in B(1)}\frac{1}{j}\log\mid f(z)\mid 
\]
holds, 
where $B(1)$ is the unit ball of ${\cal H}(j\varphi )_{C}$. 
For $r < \mbox{dist}(z,\partial\Delta^{n})$ and $f\in B(1)$, 
the mean value inequality applied to the plurisubharmonic function $\mid f\mid^{2}$ implies 
\begin{eqnarray*}
\mid f(z)\mid^{2} & \leq & \frac{1}{\pi^{n}r^{2n}/n!}\int_{\mid\zeta -z\mid < r}\mid f(\zeta )\mid^{2}d\lambda (\zeta )  \\
 &\leq & \frac{1}{\pi^{n}r^{2n}/n!}\exp (2j\sup_{\mid\zeta -z\mid < r}\varphi (\zeta ))\int_{\Delta^{n}}\mid f\mid^{2}e^{-2j\varphi}d\lambda 
\end{eqnarray*}
holds. 
If we take the supremum over all $f\in B(1)$ we have 
\[
\varphi_{j}(z) \leq \sup_{\mid\zeta -z\mid < r} \varphi (\zeta ) + \frac{1}{2j}\log \frac{1}{\pi^{n}r^{2n}/n!}
\]
holds. 

Conversely, the $L^{2}$-extension theorem (Theorem 2.6) applied twice to 
the zero dimensional subvariety $\{ z\} \subset C \subset \Delta^{n}$  shows 
that for any $a\in \mbox{\bf C}$ there is a holomorphic function 
$f$ on $\Delta^{n}$ such that 
$f(z) = a$ and 
\[
\int_{\Delta^{n}}\mid f\mid^{2}e^{-j\varphi}d\lambda 
+ 
\int_{C}\mid f\mid^{2}e^{-j\varphi}d\lambda_{C} \leq 2K_{1}\mid a\mid^{2}e^{-2j\varphi (z)},
\]
where $K_{1}$ only depends on $n$. 
We fix $a$ such that the righthandside is $1$.
This gives the other inequality 
\[
\varphi_{j}(z)\geq \frac{1}{j}\log\mid a\mid = \varphi (z) - \frac{\log K_{1}}{2j}.
\]
If $z\in \Delta^{n} - C$, 
there is a holomorphic function 
$f$ on $\Delta^{n}$ such that 
$f(z) = a$ and 
\[
\int_{\Delta^{n}}\mid f\mid^{2}e^{-j\varphi - \psi}d\lambda 
\leq K_{1}\mid a\mid^{2}e^{-2j\varphi (z)-\psi (z)}
\]
holds.
In particular $f\mid_{C}\equiv 0$ holds in this case. 
This implies the inequality
\[
\varphi_{j}(z)\geq  \varphi (z) + \frac{1}{2j}\psi (z)- \frac{\log K_{1}}{2j}.
\]
Hence we see that 
\[
\nu (\varphi_{j},z) \leq \nu (\varphi ,z)
\]
holds for every $z\in \Delta^{n}$.
In the opposite direction we find 
\[
\sup_{\mid x-z\mid <r}\varphi_{j}(x) \leq \sup_{\mid\zeta -z\mid < 2r}\varphi (\zeta ) + \frac{1}{j}\log\frac{K_{2}}{r^{n}}
\]
holds, where $K_{2}$ is a positive constant independent of $j$. 
Thus we obtain 
\[
\nu (\varphi_{j},x) \geq \nu (\varphi ,x) - \frac{n}{j}. 
\]
{\bf Q.E.D.} \vspace{5mm} \\ 

To construct $\{ h_{j}\}$ we need to globalize the above argument,i.e. 
we need to glue local approximations.  
But this is completely parallel to the argument in \cite[pp. 377-380]{d}.
Hence we omit it. 
We note that the glueing process in \cite[pp. 377-380, see especially p.377, Lemma 3.5]{d} 
does not change singularities of the sequence of approximations 
(up to quasi-isometry) on $X$ (hence in particular on $C$). 

By the construction  we have the following lemma.
\begin{lemma}
\[
{\cal I}(h^{jm})\subseteq {\cal I}(h_{j}^{jm})
\]
holds for every $j$ and $m \geq 0$.
\end{lemma}
{\bf Proof of Lemma 2.4}.
By the construction of $h_{j}$ we see that 
\[
{\cal I}(h^{j}) \subseteq {\cal I}_{\infty}(h_{j}^{j})
\]
holds (for the definition of ${\cal I}_{\infty}$ see Section 2.2). 
By the subadditivity theorem (\cite{e-d-l}),
we see that 
\[
{\cal I}(h^{jm}) \subseteq {\cal I}(h^{j})^{m}
\subseteq {\cal I}_{\infty}(h_{j}^{j})^{m} \subseteq 
{\cal I}_{\infty}(h_{j}^{jm}) \subseteq {\cal I}(h_{j}^{jm})
\]
hold for every $m\geq 0$.
\vspace{5mm} {\bf Q.E.D.} \\

By Lemma 2.4 the sequence $\{ h_{j}\}$ satisfies the 
4-th condition above. 
The 3-rd and 5-th conditions are satisfied by the convergences of the 
Lelong numbers 
\[
\lim_{j\rightarrow\infty} \nu (\Theta_{h_{j}}) 
= \nu (\Theta_{h})
\]
and 
\[
\lim_{j\rightarrow\infty} \nu (\Theta_{h_{j}}\mid_{C}) 
= \nu (\Theta_{h}\mid_{C})
\]
which follow from  Proposition 2.3. 
Since the first and the second conditions are cleary satisfied,
$\{ h_{j}\}$ is a desired sequence  of singular hermitian 
metrics on $L$. 

We note that for every $m\geq 0$, ${\cal O}_{C}(mL)\otimes{\cal I}(h^{m})$ 
is torsion free, since it is a subsheaf of a locally free
sheaf on a smooth variety $C$. 
Since $\dim C = 1$, this means that for every $m\geq 0$  
${\cal O}_{C}(mL)\otimes{\cal I}(h^{m})$ 
is invertible  on $C$. 
Since for every $0\leq k < j$ 
\[
\deg_{C}{\cal O}_{C}((jm+k) L)\otimes{\cal I}(h^{jm+k}) 
\leq 
\deg_{C}{\cal O}_{C}(jmL)\otimes{\cal I}(h^{jm}) 
+ k(L\cdot C)
\]
holds, by Lemma 2.4 we see that for every $0\leq k < j$ 
\[
\deg_{C}{\cal O}_{C}((jm+k) L)\otimes{\cal I}(h^{jm+k}) 
\leq 
\deg_{C}{\cal O}_{C}(jmL)\otimes{\cal I}(h_{j}^{jm})
+ k(L\cdot C)
\]
hold.
Then by the Riemann-Roch theorem and the Kodaira vanishing theorem 
imply that 
\[
(L,h_{j})\cdot C \geq (L,h)\cdot C
\]
holds. 
In particular we see that 
\[
\overline{\lim}_{j\rightarrow\infty}(L,h_{j})\cdot C 
\geq  (L,h)\cdot C
\]
holds.

On the other hand since $h_{j}$ has algebraic singularities, 
\[
(L,h_{j})\cdot C = (L,h_{j})*C
\]
holds.
This implies that
\[
\overline{\lim}_{j\rightarrow\infty}(L,h_{j})\cdot  C 
= \overline{\lim}_{j\rightarrow\infty}(L,h_{j})* C 
= (L,h)* C
\]
hold. 
The last equality comes from the 2-nd condition.
Combining the above inequalities, we have that 
\[
(L,h)\cdot C \leq (L,h)*C
\]
holds. 
Since we already have  the opposite inequality,  we see that 
\[
(L,h)\cdot C = (L,h)*C
\]
holds.
This completes the proof of Theorem 2.7. {\bf Q.E.D.} 
\begin{corollary}
Let $(L,h)$ be a pseudoeffective singular hermitian line bundle on 
a smooth projective variety $X$.
Let $Y$ be a subvariety such that the restriction 
$h\mid_{Y}$ is well defined. 
Then for every irreducible curve $C$ on $Y$ such that $h\mid_{C}$
is well defined, 
\[
(L,h)\cdot C = (L,h)\mid_{Y}\cdot C
\]
holds. 
In other words, the intersection theory is 
compactible with restrictions. 
In particular $(L,h)\mid_{Y}$ is numerically trivial,
if and only if $(L,h)$ is numerically trivial on $Y$. 
\end{corollary}
By the additivity of Lelong numbers we have the 
following corollary. 
\begin{corollary} 
Let $(L,h), (L^{\prime},h^{\prime})$ be singular hermitian line bundles on 
a smooth projective variety $X$ such that the 
curvature currents $\Theta_{h}$,$\Theta_{h^{\prime}}$ are positive.
Then for an irreducible curve $C$ such that 
$h\mid_{C}$ and $h^{\prime}\mid_{C}$ are both well defined, 
\[
(L\otimes L^{\prime},h\cdot h^{\prime})\cdot C = (L,h)\cdot C + (L^{\prime},h^{\prime})\cdot C
\]
holds. 
\end{corollary}
\begin{theorem} 
Let $(L,h)$ be a singular hermitian line bundle on a smooth projective 
variety $X$. 
Suppose that $\Theta_{h}$ is bounded from below by some 
negative multiple of a $C^{\infty}$-K\"{a}hler form on $X$. 
Let $D$ be a smooth divisor on $X$. 
If $h\mid_{D}$ is well defined, then 
\[
\bar{\cal I}_{D}(h) = \bar{\cal I}(h\mid_{D})
\]
holds. 
\end{theorem}
{\bf Proof of Theorem 2.8}.
Since the statement is local, we may assume that 
$X$ is the unit open polydisk $\Delta^{n} =
\{ (z_{1},\ldots ,z_{n})\in \mbox{\bf C}^{n};
\mid z_{i}\mid < 1 , 1\leq i\leq n\}$,
$D$ is the divisor $(z_{n})$ and
$L$ is a trivial bundle with singular hermitian metric 
$e^{-\varphi}$, where $\varphi$ is a plurisubharmonic 
function on $\Delta^{n}$. 

The proof of Theorem 2.8 is parallel to that of Theorem 2.7 , if we replace 
the  curve $C$ by the divisor $D$.

First we shall consider the case that $h$ has algebraic singularities. 
Let 
\[
f : Y \longrightarrow X
\]
be a modification such that there exits an effective {\bf Q}-divisor 
$F$ with normal crossings on $Y$ such that 
\[
{\cal I}(f^{*}h^{m}) = {\cal O}_{Y}(-[mF])
\]
holds for every $m\geq 0$. 
Let $E$ be the strict transform of $D$ in $Y$.
By the assumption the support of $F$ does not contain $E$. 
We may and do assume that $E+F$ is a divisor with normal crossings. 
Then we have that 
\[
{\cal I}_{E}(f^{*}h^{m}) = {\cal O}_{E}(-[mF])
\]
holds for every $m \geq 0$. 
And 
\[
(\star ) \hspace{10mm}\tilde{\cal I}_{D}(h^{m}) = {\cal O}_{X}(-K_{X})\otimes 
f_{*}({\cal O}_{Y}(K_{Y})\otimes{\cal O}_{Y}(-[mF]))\otimes{\cal O}_{D}
\]
holds.
Let 
\[
F = \sum_{i=1}^{\ell}a_{i}F_{i}
\]
be the irreducible decomposition of $F$. 
Let us fix an arbitrary point $x$ on $X$. 
For every $1\leq i\leq \ell$ and $m\geq 1$, we define 
the number $b_{i}(m)$ by 
\[
b_{i}(m): = \inf_{\sigma}\mbox{mult}_{F_{i}}f^{*}(\sigma ),
 \]
where $\sigma$ runs all the nonzero element 
of $\tilde{\cal I}_{D}(h^{m})_{x}$.

Let ${\cal H}(m\varphi )$ be the Hilbert space 
defined by 
\[
{\cal H}(m\varphi ) := \{ \phi\in {\cal O}(\Delta^{n})\mid 
\int_{\Delta^{n}}\mid \phi\mid^{2}e^{-m\varphi}d\lambda < \infty\} ,
\] 
with the inner product 
\[
(\phi ,\phi^{\prime}) := \int_{\Delta^{n}}\phi\cdot\bar{\phi}^{\prime}\,\,e^{-m\varphi}d\lambda ,
\]
where $d\lambda$ is the usual Lebesgue measure on $\Delta^{n}$.
Let $\{ \sigma_{\ell}\}$ be an orthonormal basis 
of ${\cal H}(m\varphi )$ and let 
\[
\varphi_{m}:= \frac{1}{2m}\log\sum_{\ell}\mid\sigma_{\ell}\mid^{2}.
\]
Clearly 
\[
b_{i}(m) = m\cdot \nu (f^{*}\varphi_{m},F_{i})
\]
holds.
We define the nonnegative numbers $\{ r_{i}\}$ by 
\[
K_{Y} = f^{*}K_{X} +\sum_{i}r_{i}F_{i} + \mbox{other components}.
\]
To estimate $b_{i}(m)$ we shall prove the following lemma.  
\begin{lemma}
\begin{enumerate}
\item $\nu (f^{*}\varphi_{m},F_{i}) \leq a_{i}$ holds,
\item \[
\nu (f^{*}\varphi_{m},F_{i}) \geq a_{i} - \frac{n+r_{i}}{m}
\]
holds for every $m\geq 1$. 
\end{enumerate}
\end{lemma}
{\bf Proof of Lemma 2.5.} 
The first assertion follows from the parallel argument as in 
the proof of Proposition 2.3.
In fact the $L^{2}$-extension theorem (Theorem 2.6) applied to 
the zero dimensional subvariety $\{ z\}\subset \Delta^{n}$  shows 
that for any $a\in \mbox{\bf C}$ there is a holomorphic function 
$f$ on $\Delta^{n}$ such that 
$f(z) = a$ and 
\[
\int_{\Delta^{n}}\mid f\mid^{2}e^{-j\varphi}d\lambda 
 \leq 2K_{1}\mid a\mid^{2}e^{-2j\varphi (z)},
\]
where $K_{1}$ only depends on $n$. 
This gives the inequality :
\[
\varphi_{m} \geq \frac{1}{m}\log \mid a\mid 
= \varphi -\frac{\log K_{1}}{2m}.
\]
This implies the first assertion.

Let us prove the second assertion.
Let $y\in Y$ be a general point on $F_{i}$ and let
$(U,w_{1},\ldots ,w_{n})$ be a local coordinate around $y$ such that 
$U$ is biholomorphic to the unit open polydisk in $\mbox{\bf C}^{n}$ 
by the coordinate $(w_{1},\ldots ,w_{n})$. 
We set 
\[
J = \frac{f^{*}dz_{1}\wedge \cdots \wedge dz_{n}}{dw_{1}\wedge\cdots \wedge dw_{n}}.
\]
Then $J$ is a holomorphic function on $U$.
Let $\phi \in {\cal H}(m\varphi )$ be an arbitrary 
element.

Then my mean value inequality 
\begin{eqnarray*}
\mid f^{*}\phi (w)\cdot J(w)\mid^{2}
& \leq  & \frac{1}{\pi^{n}r^{2n}/n!}\int_{\mid\zeta -w\mid <r}\mid f^{*}\phi\mid^{2}\mid J\mid^{2}d\lambda (\zeta ) \\
& \leq & \frac{1}{\pi^{n}r^{2n}/n!}\exp (2m\sup_{\mid\zeta -w\mid < r}f^{*}\varphi (\zeta ))\int_{\Delta^{n}}\mid f^{*}\phi (\xi )\cdot J\mid^{2}f^{*}e^{-2m\varphi}d\lambda (\xi ) \\
& \leq & \frac{1}{\pi^{n}r^{2n}/n!}\exp (2m\sup_{\mid\zeta -w\mid < r}f^{*}\varphi (\zeta ))\int_{\Delta^{n}}\mid\phi (z)\mid^{2}e^{-2m\varphi}d\lambda (z) 
\end{eqnarray*}
hold, where $d\lambda$ denotes the usual Lebesgue measure
on the unit open polydisk. 

We note that 
\[
\varphi_{m}(z) = \sup_{f\in B(1)}\frac{1}{m}\log\mid f(z)\mid 
\]
holds, 
where $B(1)$ is the unit ball of ${\cal H}(m\varphi )$. 
If we take the supremum over all $\phi$ in $B(1)$ in ${\cal H}(m\varphi )$, we have that 
\[
\varphi_{m}(w)  \leq \sup_{\mid\zeta -w\mid <r}\varphi (\zeta ) + \frac{1}{2m}
\log \frac{1}{\pi^{n}\mid J(w)\mid^{2}\cdot r^{2n}/n!}
\]
holds. 
Hence we have that 
\[
\nu (f^{*}\varphi_{m},y) \geq \nu (f^{*}\varphi ,y) - \frac{n+r_{i}}{m}
= a_{i} -  \frac{n+r_{i}}{m}
\]
hold. 
This completes the proof of Lemma 2.5. \vspace{5mm} {\bf Q.E.D.} \\
We note that for every positive number $\epsilon$ and positive integer $m$,
\[
{\cal I}(h^{1+2\epsilon}\mid_{D}) \subseteq \sqrt[m]{\tilde{\cal I}_{D}(h^{(1+\epsilon )m})}
\]
holds by the formula $(\star )$ and the definition of $\sqrt[m]{\tilde{\cal I}_{D}(h^{m})}$.
Hence we have that
\[
{\cal I}(h^{1+2\epsilon }\mid_{D}) \subseteq {\cal I}_{D}(h^{1+\epsilon })
\]
holds.
By the definition of the closure of multiplier ideal sheaves,
letting $\epsilon$ tend to $0$,
 we have that 
\[
\bar{\cal I}(h\mid_{D}) \subseteq \bar{\cal I}_{D}(h)
\]
holds. 

On the other hand 
by Lemma 2.5 we have that 
\[
\lim_{m\rightarrow\infty}\frac{1}{m}b_{i}(m) =
\lim_{m\rightarrow\infty}\nu (f^{*}\varphi_{m},F_{i}) = a_{i}
\]
hold. 
By the definitions of $b_{i}(m)$ and  $\sqrt[m]{\tilde{\cal I}(h^{m})}$,
we see that the opposite inclusion : 
\[
\bar{\cal I}(h\mid_{D}) \supseteq \bar{\cal I}_{D}(h)
\]
holds. 
Hence  
\[
\bar{\cal I}_{D}(h) = \bar{\cal I}(h\mid_{D})
\]
holds. 

If $h$ is not of algebraic sigularities, by approximating $h$ by a sequence 
of singular hermitian metrics with algebraic singularities as in Section 2.6,
we completes the proof of Theorem 2.8. {\bf Q.E.D.}
\section{Characterization of numerically trivial singular hermitian line bundles}
In this section we prove Theorem 1.2. 
Let $(L,h)$ be a singular hermitian line bundle 
on a smooth projective variety $X$ with positive curvature 
current. 
Suppose that $(L,h)$ is numerically trivial on $X$. 
Let us define the closed positive current $T$ on $X$ 
by 
\[
T := \frac{1}{2\pi}\Theta_{h} - \sum_{D}\nu (\Theta_{h},D)D
\]
where $D$ runs all the prime divisors on $X$. 

Let us define the subset $S$ of $X$ by  
\[
S : = \{ x\in X\mid \nu (T,x) > 0\}.
\]
Then $S$ consists of at most countable union of subvarieties 
of codimension greater than or equal to $2$ by a theorem of Siu (\cite{s}).
Let $n$ be the dimension of $X$. 
Let $H$ be a very ample divisor and let $C$ be a very general 
complete intersection curve of $(n-1)$-members of $\mid H\mid$.
If we take $H$ sufficiently ample and take $C$ very general 
we may assume that 
\[
C \cap S = \emptyset
\]
holds and $C$ intersects every prime divisor 
$D$ with $\nu (\Theta_{h},D) > 0$ (such prime divisors 
are at most coutably many) at $D_{reg}$ 
transversally. 
Let $\omega$ be a K\"{a}hler form which represents $c_{1}(H)$. 
Let $\sigma\in \Gamma (X,{\cal O}_{X}(H))$ be 
a very general nonzero element such that 
$D = (\sigma )$ is smooth and $T\mid_{D}$ is well defined. 
Then by Stokes' theorem, 
\[
T(\omega^{n-1}) = \int_{D}T\wedge\omega^{n-2}
\]
holds. 
Hence inductively we have that 
\[
T (\omega^{n-1}) = \int_{C}T 
= L\cdot C -\sum_{D}\nu(\Theta_{h},D)D\cdot C
\]
hold.
On the other hand, by the choice of $C$ and 
Lemma 2.2 we see that 
\[
{\cal I}(h^{m}\mid_{C})
\supseteq {\cal O}_{C}(-[m\sum_{D}\nu (\Theta_{h},D)D])
\]
holds for every $m\geq 0$ (since $C$ is smooth, 
the both sides are torsion free). 
Hence if 
\[
L\cdot C -\sum_{D}\nu(\Theta_{h},D)D\cdot C > 0
\]
holds, then 
\begin{eqnarray*}
(L,h)\cdot C & \geq & L\cdot C -\sum_{D}\nu(\Theta_{h},D)D\cdot C  \\
& = & \overline{\lim}_{m\rightarrow\infty}\frac{1}{m}\deg_{C}{\cal O}_{C}(mL)\otimes{\cal O}_{C}(-[m\cdot\sum \nu (\Theta_{h},D)D])  \\
& = & L\cdot C - \sum_{D}\nu (\Theta_{h},D)D\cdot C > 0
\end{eqnarray*}
hold by the Riemann-Roch theorem and the Kodaira vanishing 
theorem. 
This is the contradiction. 
Hence we see that 
\[
T (\omega^{n-1}) = \int_{C}T 
= L\cdot C -\sum_{D}\nu(\Theta_{h},D)D\cdot C = 0
\]
hold.
Since $T$ is closed positive, this implies that $T\equiv 0$.
Hence we conclude that 
\[
\Theta_{h} = 2\pi\sum_{D}\nu (\Theta_{h},D)D
\]
holds.  
This completes the proof of Theorem 1.2.  {\bf Q.E.D.} \vspace{10mm} \\
By the proof of Theorem 1.2, we obtain the following. 
\begin{theorem}
Let $(L,h)$ be a  pseudoeffective singular hermitian line bundle on a smooth projective 
variety $X$.
Then $(L,h)$ is numerically trivial if and only if for every 
irreducible curve $C$ such that the restriction $h\mid_{C}$ is 
well defined 
\[
\Theta_{h}\mid_{C} = 2\pi\sum_{x\in C}\nu (\Theta_{h}\mid_{C},x)x
\]
holds. 
\end{theorem} 
By using the intersection theory on smooth divisors (cf. Section 2.5),
we have the following corollary. 
\begin{corollary}
Let $X$ be a smooth projective variety and let 
$(L,h)$ be a pseudoeffective singular hermitian line bundle on $X$.
Let $D$ be a smooth divisor on $X$. 
Suppose that $(L,h)$ is numerically trivial on $D$. 
Then 
\[
S_{D}:= \{ x\in D\mid \nu_{D}(\Theta_{h},x) > 0\}
\]
is a sum of countably many prime divisors on $D$. 
And for every $m\geq 0$, 
\[
\bar{\cal I}_{D}(h^{m})= {\cal I}(m\cdot\sum_{E}\nu_{D}(\Theta_{h},E)E)
\] 
holds, where $E$ runs all prime divisors on $D$. 
\end{corollary}
Here we do not need assume that the restriction  $h\mid_{D}$ 
is well defined.  \vspace{5mm} \\
{\bf Proof of Corollary 3.1}.
The proof is essentially same as that of Theorem 1.2. 

Let $X$,$D$,$(L,h)$ be as above. 
Let $\{ F_{i}\}_{i\in I}$ be the set of divisorial components of 
$S_{D}$. 
Let $C$ be a {\bf very general} complete intersection curve 
of a sufficiently ample linear system $\mid H\mid$ on $D$ which 
does not intersects $S_{D} - \cup_{i\in I}F_{i}$ and 
meets every $F_{i} (i\in I)$ transversally.  
We set 
\[
a_{i} := \nu_{D}(\Theta_{h}.F_{i}).
\]
By the definition of the intersection number
\[
(L,h)\cdot C = (L -\nu (\Theta_{h},D)D)\cdot C
+\overline{\lim}_{m\rightarrow\infty}m^{-1}\deg_{C}\tilde{\cal I}_{D}(h^{m})
\otimes{\cal O}_{C} 
= 0
\]
hold.
Hence  if we take $C$ very general, we see that 
\[
(L -\nu (\Theta_{h}.D)D)\cdot C = \sum_{i\in I}a_{i}(F_{i}\cdot C)
\]
holds by the definition of $\nu_{D}$. 

Suppose that $S_{D} - \cup_{i\in I}F_{i}$ is nonempty. 
Let $\{ C_{t}\}_{t\in\Delta}$ be a family of complete intersection 
curve of $(\dim D - 1)$-members of $\mid H\mid$
such that 
\begin{enumerate}
\item $C_{0}$ is not contained in $S_{D}$,
\item $C_{0}\cap (S_{D} -\cup_{i\in I}F_{i})\neq \emptyset$,
a very general member of $\{ C_{t}\}$  does not intersects $S_{D} - \cup_{i\in I}F_{i}$ and 
meets every $F_{i}$ transversally.
\end{enumerate}
Then by the uppersemicontinuity of Lelong numbers $\nu_{D}$ 
in countable Zariski topology (the uppersemicontinuity is ovbious 
by the definition of $\nu_{D}$ (cf. Section 2.5)), we see that 
\[
(L,h)\cdot C_{0} < 0
\]
holds. 
This is the contradiction, since $(L,h)\cdot C_{0}\geq 0$ holds by the 
definition of the intersection number.  
Hence we have that $S_{D} = \sum_{i\in I} F_{i}$ holds. 
This completes the proof of the first assertion. 

Let us prove the second assertion.
Let us fix an irreducible component $F_{i}$ of $S_{D}$. 
Let $U$ be a Stein open subset on $D$ and let $\sigma \in \tilde{\cal I}(h^{\ell})(U)$ be a very general element.
Let $\epsilon$ be any small positive number. 
Then by the definition of $\nu_{D}(\Theta_{h})$, we see that
for every sufficiently large $\ell$ 
\[
\mbox{($\star\star$)}\hspace{10mm}(1-\epsilon )\nu_{D}(\Theta_{h},F_{i})
\leq \frac{1}{\ell}\mbox{mult}_{F_{i}}(\sigma) 
\leq  \nu_{D}(\Theta_{h},F_{i})
\]
hold. 
Then by the definition of ${\cal I}_{D}(h^{m})$ and Lemma 2.2, we see that 
\[
\bar{\cal I}_{D}(h^{m})\subseteq {\cal I}(m\cdot\sum_{i\in I}\nu_{D}(\Theta_{h},F_{i})F_{i})
\] 
holds.

Let $C$ be a very general smooth complete intersection of 
$(\dim D-1)$-members of $\mid H\mid$. 
Then as above
\[
((L-\nu (\Theta_{h},D)D)\mid_{D} - \sum_{i\in I}a_{i}F_{i})\cdot C
= 0
\]  
holds.

We claim that 
$L - \nu (\Theta_{h},D)D -\sum_{i\in I}a_{i}F_{i}$ 
is pseudoeffective in the sense that 
$c_{1}((L - \nu (\Theta_{h},D)D)\mid_{D} -\sum_{i\in I}a_{i}F_{i})$
is on the closure of effective cone of $D$. 
Let $G$ be an ample line bundle on $X$ 
such that 
${\cal O}_{X}(G +mL)\otimes {\cal I}(h^{m})$ is 
globally generated on $X$ for every $m\geq 0$. 
This is possible by \cite[p.664, Proposition 1]{si}. 
By the formula ($\star\star$), this implies that for any sufficiently small 
$\epsilon > 0$ and every finite subset $I_{0}$ of $I$
\[
(L - \nu(\Theta_{h},D)D)\mid_{D}
- \sum_{i\in I_{0}}a_{i}F_{i}
\]
is pseudoeffective. 
Hence this we see that $(L - \nu (\Theta_{h},D)D)\mid_{D} -\sum_{i\in I}a_{i}F_{i}$
is pseudoeffective. 

Since 
\[
((L - \nu (\Theta_{h},D)D)\mid_{D} -\sum_{i\in I}a_{i}F_{i})
\cdot H^{\dim D-1} = 0
\]
holds, 
this implies that $L - \nu (\Theta_{h},D)D)\mid_{D} -\sum_{i\in I}a_{i}F_{i}$
is numerically trivial.

Let $f : \tilde{D} \longrightarrow D$ 
be any composition of successive blowing ups with smooth 
center, 
then by the same argument as above, we see that 
\[
f^{*}(L - \nu (\Theta_{h},D)D)\mid_{D}
- \sum_{\tilde{E}}\nu_{\tilde{D}}(f^{*}\Theta_{h},\tilde{E})\tilde{E}
\]
is numerically trivial, where $\tilde{E}$ runs all the 
prime divisors on $\tilde{D}$. 
We note that by the definitions of $\nu_{D}$ and $\nu_{\tilde{D}}$,
we see that 
\[
\sum_{\tilde{E}}\nu_{\tilde{D}}(f^{*}\Theta_{h},\tilde{E})\tilde{E}
- f^{*}(\sum_{i\in I}a_{i}F_{i})
\]
is effective, i.e. a sum of prime divisors with nonnegative
coefficients. 
Since $f^{*}((L - \nu (\Theta_{h},D)D)\mid_{D} -\sum_{i\in I}a_{i}F_{i})$
is numerically trivial on $\tilde{D}$,
we see that 
\[
\sum_{\tilde{E}}\nu_{\tilde{D}}(f^{*}\Theta_{h},\tilde{E})\tilde{E}
 = f^{*}(\sum_{i\in I}a_{i}F_{i}). 
\]
Let $m$ be any positive integer and  $f_{m} : D_{m}\longrightarrow D$ be a modification such that 
$f^{*}_{m}\bar{\cal I}(h^{m})$  is locally free. 
Then by the definition of $\bar{\cal I}_{D}(h^{m})$, 
it is determined by the Lelong numbers $\nu_{D_{m}}$ 
on prime divisors on $D_{m}$. 
Applying the above argument by taking $\tilde{D}$ to 
be $D_{m}$, we see that 
\[
\bar{\cal I}_{D}(h^{m}) = {\cal I}(m\cdot\sum_{i\in I}a_{i}F_{i})
\]
holds.
This completes the proof of Corollary 3.1.
{\bf Q.E.D.}
\begin{remark}
By the above proof Corollary 3.1 still holds for a subvariety $V$ on 
$D$, if there exists a curve on $V$ such that 
$(L,h)\cdot C$ is well defined (cf. \cite[Remark 3.1]{tu3}). 
\end{remark}
  
\section{Numerical triviality and the growth of $H^{0}$}
In this section we shall relate the numerical triviality of 
singular hermitian line bundles with positive curvature current 
and the growth of dimension of global sections. 
\begin{definition}
Let $(L,h)$ be a singular hermitian line bundle on a 
smooth projective variety $X$.
Let $H$ be an ample  line bundle on $X$.
We define the number $\mu_{h}(X,H+mL)$ by 
\[
\mu_{h}(X,H+mL) := (\dim X)!\cdot\overline{\lim}_{\ell\rightarrow\infty}\ell^{-\dim X}
\dim H^{0}(X,{\cal O}_{X}(\ell (H+mL))\otimes{\cal I}(h^{m\ell}))
\]
For a subvariety $Y$ in $X$ such that $(L,h)\mid_{Y}$ is well defined, we define 
\[
\mu_{h}(Y,H+mL) := (\dim Y)!\cdot\overline{\lim}_{\ell\rightarrow\infty}\ell^{-\dim Y}
\dim H^{0}(Y,{\cal O}_{Y}(\ell (H+mL))\otimes{\cal I}(h^{m\ell})/tor),
\]
where $tor$ denotes the torsion part of ${\cal O}_{Y}(H+mL)\otimes{\cal I}(h^{m})$. 
\end{definition}
We note that $\mu_{h}(Y,H+mL)$ is different from 
$\mu_{h}(Y,H+mL\mid_{Y})$ in general. 
By Corollary 2.1, we note that 
if $(L,h)$ is numerically trivial on $Y$ if and only if  
$(L,h)\mid_{Y}$ is numerically trivial. 

\begin{lemma}
Suppose that $(L,h)$ is pseudoeffective and is  not numerically trivial on $X$. 
Then 
\[
\overline{\lim}_{m\rightarrow\infty}m^{-1}\mu_{h} (X,H+mL) > 0 
\]
holds for every ample line bundle $H$ on $X$. 
\end{lemma}
{\bf Proof of Lemma 4.1.}
Let $n$ be the dimension of $X$.
We prove this lemma by induction on $n$.
If $n =1$, then  for every $\ell\geq 0$
\[
{\cal O}_{X}(\ell L)\otimes{\cal I}(h^{\ell})
\supseteq 
{\cal O}_{X}(\ell L -\sum_{x\in X}[\ell\cdot\nu (\Theta_{h},x)])
\]
holds by Lemma 2.2. 
Hence by Theorem 1.2, we see that 
\[
\lim_{\ell\rightarrow\infty}
\deg_{X}{\cal O}_{X}(\ell L)\otimes{\cal I}(h^{\ell}) 
= +\infty .
\]
This implies Lemma 4.1. 

Let $\pi :\tilde{X} \longrightarrow {\bf P}^{1}$
be a Lefschetz pencil associated with a very ample linear 
system  say $\mid H\mid$ on $X$.
If we take the pencil very general, we may assume that 
${\cal I}(h^{\ell})$ is an ideal sheaf on 
all fibers of $\pi$ for every $\ell$.
Let 
\[ 
b : \tilde{X}\longrightarrow X
\]
be the modification associated with the pencil and let
$E$ be the exceptional locus of $b$.
We note that on the Hilbert scheme of curves in $X$, 
the intersection number $(L,h)\cdot C$ is lower semicontinuous 
in countable Zariski topology
by the upper semicontinuity of the Lelong number (or by the $L^{2}$-extension theorem (Theorem 2.6)), where $C$ moves in the Hilbert scheme. 
Then by the inductive assumption for a general fiber $F$ of
$\pi$ we see that 
\[
\overline{\lim}_{m\rightarrow\infty}m^{-1}\mu_{h} (F,b^{*}(H+mL)) > 0
\]
holds.
Let us consider the direct image
\[
{\cal E}_{m,\ell} := \pi_{*}{\cal O}_{\tilde{X}}(\ell b^{*}(H+mL)\otimes
{\cal I}(b^{*}(h^{m\ell}))).
\]
By Grothendiek's theorem, we see that
\[
{\cal E}_{m,\ell} \simeq \oplus_{i=1}^{r}{\cal O}_{{\bf P}^{1}}(a_{i})
\]
for some $a_{i} = a_{i}(m,\ell )$ and $r = r(m,\ell )$.
By the inductive assumption we see that 
\[
\overline{\lim}_{m\rightarrow\infty}m^{-1}(\overline{\lim}_{\ell\rightarrow\infty}
\ell^{-(n-1)}r(m,\ell)) > 0
\]
holds.
We note that $\ell_{0}b^{*}H -E$ is ample
 for some large positive integer $\ell_{0}$.
Hence we see that 
\[
{\cal O}_{\tilde{X}}(\ell_{0}b^{*}H - E)
\]
admits a $C^{\infty}$-hermitian metric $h_{0}$ with strictly positive curvature.Let $h_{1}$ be a $C^{\infty}$-hermitian metric 
on ${\cal O}_{\mbox{{\bf P}}^{1}}(1)$. 
Then there exists a positive rational number $c$ such that 
\[
\frac{1}{\ell_{0}}\Theta_{h_{0}}-
c\cdot\pi^{*}\Theta_{h_{1}}
\]
is a K\"{a}hler form on $\tilde{X}$. 
By Nadel's vanishing theorem (Theorem 2.1),
\[
H^{1}(\tilde{X},{\cal O}_{\tilde{X}}(\ell (b^{*}(H+mL)-\frac{1}{\ell_{0}}E))
\otimes{\cal I}(b^{*}(h^{m\ell}))
\otimes \pi^{*}{\cal O}_{{\bf P}^{1}}(- c\ell )) = 0
\]
holds for every sufficiently large $\ell$  such that $\ell /\ell_{0}$
and $c\ell$ are integers. 
Also by Nadel's vanishing theorem, we see that 
\[
R^{1}\pi_{*}{\cal O}_{\tilde{X}}(\ell (b^{*}(H+mL)-\frac{1}{\ell_{0}}E))
\otimes{\cal I}(b^{*}(h^{m\ell}))
\]
is the $0$-sheaf on $\mbox{\bf P}^{1}$ for every 
sufficiently large $\ell$ divisible by $\ell_{0}$.
Hence we see that 
${\cal E}_{m,\ell}\otimes {\cal O}_{\mbox{\bf P}^{1}}(-c\ell +1)$ is  
globally generated on $\mbox{\bf P}^{1}$ 
for every sufficiently large $\ell$  such that $\ell /\ell_{0}$
and $c\ell$ are integers. 
This implies that 
\[
\overline{\lim}_{\ell\rightarrow\infty}\ell^{-1}\min_{i}a_{i} \geq c
\]
holds for every $i$.
Hence 
\[
\overline{\lim}_{\ell\rightarrow\infty}
\ell^{-n}\dim H^{0}(\tilde{X},{\cal O}_{\tilde{X}}(\ell b^{*}(H+mL))
\otimes{\cal I}(b^{*}(h^{m\ell}))) \geq 
\]
\[
\,\,\,\,\,\,\,\,\,\,\,\,\, c\cdot\overline{\lim}_{\ell\rightarrow\infty}
\ell^{-(n-1)}r(m,\ell)
\] 
holds.
By this we  see that 
\[
\overline{\lim}_{m\rightarrow\infty}m^{-1}(\overline{\lim}_{\ell\rightarrow\infty}
\ell^{-n}\dim H^{0}(\tilde{X},{\cal O}_{\tilde{X}}(\ell b^{*}(H+mL))
\otimes{\cal I}(b^{*}(h^{m\ell}))) > 0
\]
holds.
Since 
\[
b_{*}{\cal I}(b^{*}h^{m\ell}) \subseteq
{\cal I}((h^{m\ell}))
\]
holds by Lemma 2.1,
we see that
\[
\overline{\lim}_{m\rightarrow\infty}m^{-1}\mu_{h} (X,H+mL) > 0
\]
holds.
Here we have assumed that $H$ to be sufficiently very ample. 
To prove the general case of Lemma 4.1, we argue as follows.
Let $H$ be any ample line bundle on $X$.
Then  thanks to Nadel's vanishing theorem
\[
\mu_{h}(X,a(H+mL)) = a^{n}\cdot \mu_{h}(X,H+mL)
\]
holds for every positive integer $a$. 
Now it is clear that Lemma 4.1 holds for any ample line bundle $H$. 
This completes the proof of Lemma 4.1. {\bf Q.E.D.}
\begin{theorem}
Let $(L,h)$ be a pseudoeffective singular hermitian line bundle on 
a smooth projective variety $X$. 
Then $(L,h)$ is numerically trivial if and only if
\[
\overline{\lim}_{m\rightarrow\infty}\mu_{h}(X,H+mL) < \infty
\]
holds for every ample line bundle $H$ on $X$. 
\end{theorem} 
{\bf  Proof of Theorem 4.1}. 
By Lemma 4.1, $(L,h)$ is numerically trivial, if 
\[
\overline{\lim}_{m\rightarrow\infty}\mu_{h}(X,H+mL) < \infty
\]
holds for every ample line bundle $H$ on $X$.  

Let us prove the converse. 
Suppose that 
\[
\overline{\lim}_{m\rightarrow\infty}\mu_{h}(X,H+mL) =  \infty
\]
holds for some ample line bundle $H$ on $X$. 
Let $x$ be a very general point of $X$ such that 
\[
{\cal I}(h^{m})_{x} = {\cal O}_{X,x}
\]
holds for every $m\geq 0$. 
\begin{lemma}
For every postive integer $N$ there exists a positive integer $m_{0}$ such that for every sufficienly large $\ell$ there exists a section 
\[
\sigma_{\ell} \in H^{0}(X,{\cal O}_{X}(\ell (H+m_{0}L))\otimes {\cal I}(h^{m_{0}\ell})\otimes {\cal M}_{x}^{N\ell }) - \{ 0\}.
\]
\end{lemma}
{\bf Proof of Lemma 4.2.} 
\[
\overline{\lim}_{m\rightarrow\infty}\mu_{h}(X,H+mL) = \infty
\]
holds by the assumption.
Hence there exists a positive integer $m_{0}$ such that 
\[
\mu_{h} (X,H+m_{0}L)  > N^{\dim X} +1 
\]
holds. 
Then 
\[
\dim H^{0}(X,{\cal O}_{X}(\ell (H+m_{0}L))
\otimes {\cal I}(h^{m_{0}\ell})) \geq  \frac{N^{\dim X}+1}{(\dim X)!}\ell^{\dim X} + o(\ell^{\dim X})
\]
holds. 
We consider the exact sequence 
\[
0\rightarrow   H^{0}(X,{\cal O}_{X}(\ell (H+m_{0}L))
\otimes {\cal I}(h^{m_{0}\ell})\otimes {\cal M}_{x}^{N\ell}) 
\rightarrow 
 H^{0}(X,{\cal O}_{X}(\ell (H+m_{0}L))
\otimes {\cal I}(h^{m_{0}\ell}))
\]
\[
\rightarrow 
 H^{0}(X,{\cal O}_{X}(\ell (H+m_{0}L))
\otimes {\cal I}(h^{m_{0}\ell})\otimes{\cal O}_{X}/{\cal M}_{x}^{N\ell}). 
\]
Since 
\[
{\cal I}(h^{m})_{x} = {\cal O}_{X,x}
\]
holds for every $m \geq 0$, we see that 
\[
\dim H^{0}(X,{\cal O}_{X}(\ell (H+m_{0}L))
\otimes {\cal I}(h^{m_{0}\ell})\otimes{\cal O}_{X}/{\cal M}_{x}^{N\ell}) 
= \frac{N^{\dim X}}{(\dim X)!}\ell^{\dim X} + o(\ell^{\dim X})
\]
holds.
Combining the above facts, we see that 
\[
H^{0}(X,{\cal O}_{X}(\ell (H+m_{0}L))
\otimes {\cal I}(h^{m_{0}\ell})\otimes {\cal M}_{x}^{N\ell})\neq 0 
\]
holds for every sufficiently large $\ell$. 
This completes the proof of Lemma 4.2. {\bf Q.E.D.} \vspace{5mm} \\

Let us continue the proof of Theorem 4.1.
Let $H_{0}$ be a sufficiently ample line bundle.
Let $C$ be a very general complete intersection of $(\dim X - 1)$-members of 
$\mid H_{0}\mid$ such that $x\in C$.
We may assume that $h\mid_{C}$ is well defined and 
\[
\sigma_{\ell}\mid_{C}\not{\equiv} 0
\]
holds for every sufficiently large $\ell$. 
This implies by a degree argument that 
\[
\overline{\lim}_{m\rightarrow\infty}\mu_{h}(C,H+mL) \geq N
\]
holds. 
Since $N$ is arbitrary, we may take $m_{0}$ so that
\[
\mu_{h}(C,H+m_{0}L) \geq 3H\cdot C
\]
holds. 
Then for every sufficiently large $\ell$
\[
\dim H^{0}(C,{\cal O}_{C}(\ell (H+m_{0}L)-\ell H))
\otimes {\cal I}(h^{m}))
\geq (H\cdot C)\cdot \ell 
\]
holds. 
Hence 
\[
(L,h)\cdot C > 0
\]
holds. 
This completes the proof of Theorem 4.1. {\bf  Q.E.D.} 
\begin{theorem}
Let $f : Y\longrightarrow X$ be a surjective morphism between 
smooth projective varieties. 
Let $(L,h)$ be a pseudoeffective singular hermitian line bundle on $X$.  
Then $(L,h)$ is numerically trivial on $X$ if and only if
$f^{*}(L,h)$ is numerically trivial on $Y$. 
\end{theorem}
{\bf Proof of Theorem 4.2}. If $(L,h)$ is numerically trivial on $X$, 
then by Theorem 1.2, $\Theta_{h}$ is a sum of at most 
countably many prime divisors with nonnegative 
coefficients. 
Hence $f^{*}\Theta_{h}$ is at most countably many prime 
divisors with nonnegative coefficients. 
Hence by Theorem 3.1, $f^{*}(L,h)$ is numerically trivial
on $Y$.

Suppose that $(L,h)$ is not numerically trivial on 
$X$. 
Let $H$ be a sufficiently very ample line bundle on $Y$
and let $C$ be a very general complete intersection curve
of $\dim Y -1$ members of $\mid H\mid$.
We may assume that 
\begin{enumerate}
\item $C$ is smooth,
\item $f(C)$ is a smooth curve,
\item $f\mid_{C} : C\longrightarrow f(C)$ is unramified 
on $\{ y\in C\mid \nu (f^{*}\Theta_{h}\mid_{C},y) > 0\}$,
\item $(L,h)\cdot f(C) > 0$ holds. 
\end{enumerate}
Then we have that 
\[
\frac{1}{2\pi}\int_{C}f^{*}\Theta_{h} - \sum_{y\in C}\nu (f^{*}\Theta_{h}, y)
= \deg (f\mid_{C})\cdot (\frac{1}{2\pi}\int_{f(C)}\Theta_{h} - \sum_{x\in f(C)}\nu (\Theta_{h}\mid_{C},x)) > 0
\]
holds.
Hence $f^{*}(L,h)$ is not numerically trivial on $Y$.
This completes the proof of Theorem 4.2. 
{\bf Q.E.D.} \vspace {10mm}\\
By Theorem 4.2, we may define the numerical triviality of pseudoeffective singular hermitian line bundles  on singular varieties. 
\begin{definition} 
Let $X$ be a singular variety and let 
\[
\pi : \tilde{X}\longrightarrow X
\]
be a resolution of singularities. 
Let $L$ be a line bundle on $X$. 
A hermitian metric $h$ on $L\mid_{X_{reg}}$ is said to be a singular hermitian 
metric on $X$, if $\pi^{*}h$ is a singular hermitian metric with  curvature current bounded from below by a $C^{\infty}$-from on $\tilde{X}$. 

$(L,h)$ is said to be pseudoeffective, if $\pi^{*}(L,h)$ is 
pseudoeffective.

Suppose that  $X$ is proper and $\tilde{X}$ is smooth projective . 
A singular hermitian line bundle $(L,h)$ is said to be numerically trivial, 
if $\pi^{*}(L,h)$ is numerically trivial. 
\end{definition}
The above definition is independent of the choice of the resolution $\pi$, 
by the $L^{1}$-property of almost plurisubharmonic functions 
and Theorem 4.2. 
\section{The fibration theorem} 
In this section we shall prove Theorem 1.1. 
\subsection{Key lemma}
The following lemma is the key for the proof of Theorem 1.1. 
\begin{lemma}
Let  $f : M \longrightarrow B$ be an algebraic fiber space 
and let $(L,h)$ be a  pseudoeffective singular hermitian line bundle on $M$. 
Suppose that for every very general fiber $F$, 
$(L,h)$ is numerically trivial on $F$ 
and there exists a subvariety $W$ of $M$ such that 
\begin{enumerate}
\item $h\mid_{W}$ is well defined,
\item $(L,h)$ is numerically trivial on $W$. 
\item $f(W) = B$.
\end{enumerate}
Then $(L,h)$ is numerically trivial on $M$. 
\end{lemma} 
{\bf Proof of Lemma 5.1}. 
Taking a suitable modification of $M$, 
by Theorem 1.2 and Theorem 4.2 we may assume that $W$ is a smooth divisor. 

Suppose that $(L,h)$ is not numerically trivial on $M$.
Then there exists an ample line bundle $H$ on $M$  such that 
\[
\overline{\lim}_{m\rightarrow \infty}\mu_{h}(M,mL+H) = \infty
\]
holds. 
We may assume that $H$ is very ample on $M$.
 
By the assumption we see that 
\[
{\cal I}(h^{m})_{x} = {\cal O}_{M,x}
\]
for a very general point $x \in W$ and every $m \geq 0$. 
Let $x_{0}$ be a very general point of $W$ such that 
\[
{\cal I}(h^{m})_{x_{0}} = {\cal O}_{M,x_{0}}
\]
holds for every $m \geq 0$.  
The proof of the following lemma is identical to that of 
Lemma 4.2.  Hence we omit it.
\begin{lemma}
For any positive integer $N$ 
there exists a  positive integer $m_{0}$  such that 
\[
H^{0}(M,{\cal O}_{M}(\ell (m_{0}L+H))\otimes {\cal I}(h^{\ell m_{0}})\otimes {\cal M}_{x_{0}}^{N\ell}) \neq 0
\]
holds for every sufficiently large $\ell$. 
\end{lemma}

Let us continue the proof of Lemma 5.1.
Let $N$ be a sufficiently large positive integer and let 
$m_{0}$ be the integer as in Lemma 5.2.  
For every sufficiently large $\ell$, we take an element 
\[
\sigma_{\ell}\in H^{0}(M,{\cal O}_{M}(\ell (m_{0}L+H))\otimes {\cal I}(h^{\ell m_{0}})\otimes {\cal M}_{x_{0}}^{N\ell}) - \{ 0\} . 
\]
Let ${\cal R}$ be the family of smooth curves 
which are complete intersection  of $\dim W -1$ members 
of $\mid H\mid_{W}$ on $W$. 
Let $d_{0}$ be a large positive integer such that 
\[
H^{\dim F-1}\cdot F\cdot (H-d_{0}W) < 0  
\]
holds for every general fiber $F$ of $f$.
Since $(L,h)$ is numerically trivial on $W$, we have 
the following lemma.
\begin{lemma} 
There exists a positive constant $A_{0}$ independent of $m_{0}$
such that  
for every member $R$ of ${\cal R}$ such that the restriction 
$h\mid_{R}$ is well defined, 
\[
(*) \,\,\, \dim H^{0}(R,{\cal O}_{R}(\ell (H+m_{0}L)-sW)\otimes 
{\cal I}(h^{m_{0}\ell}))
\leq A_{0}\cdot \ell + o(\ell )   
\] 
holds for every $0\leq s \leq d_{0}\ell$.
\end{lemma}
{\bf Proof of Lemma 5.3}. 
We note that ${\cal I}(h^{m_{0}\ell})\mid_{R}$ is torsion free, hence 
locally free (note that $\dim R = 1$),
since $R$ is smooth and ${\cal I}(h^{m_{0}\ell})\mid_{R}$ 
is a subsheaf of ${\cal O}_{R}$.
Then since $(L,h)\mid_{W}$ is numerically trivial, by Corollary 2.1,
by the formula $(\flat )$ in Remark 2.5 there exists a 
positive constant $A_{0}$ independent of $m_{0}$ such that 
\[
\deg_{R}({\cal O}_{R}(\ell (H+m_{0}L)-sW)\otimes 
{\cal I}(h^{m_{0}\ell})) \leq A_{0}\cdot\ell + o(\ell )
\]
holds for every $\ell \geq 0$ and $0\leq s \leq d_{0}\ell$.

First let us consider the case that $W\cdot R \leq 0$ holds. 
Since $H$ is ample, we have that 
\[
H^{1}(R, {\cal O}_{R}(\ell (H+m_{0}L)-sW)\otimes 
{\cal I}(h^{m_{0}\ell})) 
= 0
\]
holds for every sufficiently large $\ell$ and $0\leq s \leq d_{0}\ell$. 
By the Riemann-Roch theorem 
we have that 
\[
\dim H^{0}(R,{\cal O}_{R}(\ell (H+m_{0}L)-sW)\otimes 
{\cal I}(h^{m_{0}\ell}))
\leq A_{0}\cdot \ell + o(\ell )   
\] 
holds for every $0\leq s \leq d_{0}\ell$.

Next let us consider the case that $W\cdot R > 0$ holds. 
Then there exists a positive integer $a_{0}$ such that 
for every $s \geq a_{0}$, 
$H^{0}(R,{\cal O}_{R}(sW)) \neq 0$ holds. 
This implies that 
\[
\mbox{(1)} \hspace{5mm}\dim H^{0}(R,{\cal O}_{R}(\ell (H+m_{0}L)-sW)\otimes 
{\cal I}(h^{m_{0}\ell})) \leq 
\dim H^{0}(R,{\cal O}_{R}(\ell (H+m_{0}L))\otimes 
{\cal I}(h^{m_{0}\ell}))
\]
holds for every $s \geq a_{0}$.
Hence as before we see that there exists a positive constant $A_{0}$ such that 
\[
\dim H^{0}(R,{\cal O}_{R}(\ell (H+m_{0}L)-sW)\otimes 
{\cal I}(h^{m_{0}\ell})) \leq  A_{0}\cdot\ell +o(\ell)
\]
holds for every $s \geq a_{0}$. 
On the other hand since $H$ is ample, 
for every sufficiently large $\ell$ 
and every $0\leq s \leq a_{0}$,
 we see that 
\[
H^{1}(R,{\cal O}_{R}(\ell (H+m_{0}L)-sW)\otimes 
{\cal I}(h^{m_{0}\ell})) = 0
\]
holds.
Hence we see that by the Riemann-Roch theorem 
\[
\mbox{(2)} \hspace{5mm} \dim H^{0}(R,{\cal O}_{R}(\ell (H+m_{0}L)-sW)\otimes 
{\cal I}(h^{m_{0}\ell})) = 
\]
\[ 
\hspace{20mm}  = 1 - g(R) + \deg_{R}{\cal O}_{R}(\ell (H+m_{0}L)-sW)\otimes 
{\cal I}(h^{m_{0}\ell}) 
\]
\[
\hspace{15mm} \leq 1 - g(R) + \deg_{R}{\cal O}_{R}(\ell (H+m_{0}L))\otimes 
{\cal I}(h^{m_{0}\ell})
\]
hold for every sufficiently large $\ell$ and every 
$0\leq s\leq a_{0}$, where $g(R)$ denotes the genus of $R$. 
Combining the inequalities 
(1) and (2) above,  we see that there exists a positive constant $A_{0}$ 
such that 
\[
\dim H^{0}(R,{\cal O}_{R}(\ell (H+m_{0}L)-sW)\otimes 
{\cal I}(h^{m_{0}\ell}))  \leq A_{0}\cdot\ell + o(\ell)
\]
holds for every sufficiently large $\ell$ 
and $0\leq s \leq d_{0}\ell$, also in this case.  
{\bf Q.E.D.} \vspace{5mm} \\ 

Take $N > A_{0}$ and the 
corresponding $m_{0}$ in Lemma 5.2.
We see that using the case $s = 0$ of $(*)$,  
for every member $R$ of ${\cal R}$ containing $x_{0}$,
by a degree argument
\[
\sigma_{\ell} \mid_{R} \equiv 0
\]
holds for every sufficiently large $\ell$. 
Since the members of ${\cal R}$ containing $x_{0}$ 
dominates $W$, we see that 
\[
\sigma_{\ell}\mid_{W} \equiv 0
\]
holds for every sufficiently large $\ell$. 
Next we consider the vanishing order of $\sigma_{\ell}$ 
along the divisor $W$. 
Repeating the same arugument we see that 
\[
\sigma_{\ell} \in H^{0}(M,{\cal O}_{M}(\ell(H+m_{0}L)-d_{0}\ell W)\otimes 
{\cal I}(h^{m_{0}\ell}))
\]
holds. 
Let $F$ be a very general fiber of $f$ such that  $(L,h)\mid_{F}$ well defined and is numerically trivial.

Let ${\cal S}_{F}$ denote the family of 
smooth curves complete intersection
of $\dim F -1$ members of the very ample linear system $\mid H\mid_{F}$ 
on $F$. 
We note that $F$ is dominated by a family ${\cal S}_{F}$ of smooth curves passing through 
$W\cap F$ and 
\[
H^{\dim F-1}\cdot F \cdot (H-d_{0}W) < 0
\]
holds.
Since $\sigma_{\ell}\mid_{F}$ has the vanishing order at least $d_{0}\ell$ 
along $W\cap F$ and $(L,h)$ is numerically trivial on $F$, for every sufficiently large $\ell$,
$\sigma_{\ell}$ is identically $0$ along any members of ${\cal S}_{F}$.
In fact for every $[S]\in{\cal S}_{F}$ and every sufficiently large $\ell$  
\[
\deg_{S}{\cal O}_{S}(\ell(H+m_{0}L)-d_{0}\ell W)\otimes 
{\cal I}(h^{m_{0}\ell})
\]
is negative, since 
\[
\overline{\lim}_{m\rightarrow\infty}\frac{1}{m}\deg_{S}{\cal O}_{S}(mL)\otimes
{\cal I}(h^{m}) = (L,h)\cdot S = 0
\]
hold (cf. the formula $(\flat )$ in Remark 2.5). 
Hence 
\[
\sigma_{\ell}\mid_{F} \equiv 0
\]
holds for every sufficiently large $\ell$.
Moving smooth fibers $F$, ${\cal S}_{F}$ forms 
a dominating family of curves ${\cal S}$ on $M$. 
We may take such $\ell$ independent of a very general $F$, since
there exists a nonempty Zariski open subset ${\cal S}_{\ell}$ of 
${\cal S}$ such that 
for every $[S] \in {\cal S}_{\ell}$,  
${\cal I}(h^{m_{0}\ell})\otimes {\cal O}_{S}$ is 
an ideal sheaf on $S$ and 
\[
\deg_{S}{\cal O}_{S}(\ell(H+m_{0}L)-d_{0}\ell W)\otimes 
{\cal I}(h^{m_{0}\ell})
\]
is independent of $[S]\in {\cal S}_{\ell}$. 
This implies that 
\[
\sigma_{\ell} \equiv 0
\]
holds on $M$ for every sufficiently large $\ell$. 
This is the contradiction.  
This completes the proof of Lemma 5.1. 
\vspace{10mm} {\bf Q.E.D.} \\
\subsection{Proof of Theorem 1.1}
Let $x$ be an arbitrary point on $X$. 
We set 
\[
{\cal N}(x) := \{ V \mid \mbox{a subvariety of $X$ such that 
$x \in V$ and $(L,h)$ is } 
\]
\[
\,\,\,\,\,\,\,\,\,\,\,\,\,\,\,\,\,\,\,\,\,\,\,\,\,\,\,\,\,\,\mbox{numerically trivial on $V$} \}.
\] 
Let $\nu(x)$ denote the maximal dimension of the member of 
${\cal N}(x)$ and we set 
\[
\nu := \inf_{x\in X} \nu (x). 
\]
Then for very general $x\in X$, we see that $\nu = \nu (x)$ holds.
We note that on the Hilbert scheme of curves in $X$, 
the intersection number $(L,h)\cdot C$ is lower semicontinuous 
in countable Zariski topology
by the upper-semicontinuity of the Lelong number (or by the $L^{2}$-extension theorem (Theorem 2.6)), where $C$ moves in the Hilbert scheme.
Hence for every irreducible component of the Hilbert scheme of $X$, 
the set of  members on which the restriction of  $(L,h)$ is well defined and 
numerically trivial is locally closed in countable Zariski topology.  
Since the Hilbert scheme of $X$ has only countably many components, 
this implies  that there exists an irreducible subvariety  ${\cal N}^{0}$
 in the Hilbert scheme of 
$X$ whose members dominate $X$ and for a very general point $x$,
there exists a member $V$ of ${\cal N}^{0}$ such that 
\begin{enumerate}
\item $x \in V$,
\item $\dim V = \nu$,
\item $(L,h)$ is numerically trivial on $V$.
\end{enumerate} 
Let 
\[
\varphi : {\cal V} \longrightarrow {\cal N}^{0}
\]
be the universal family and let 
\[
p : {\cal V} \longrightarrow X
\]
be the natural morphism. 
\begin{lemma}
Let  $V$ be a very general member of ${\cal N}_{0}$.
Then there exists a Zariski open subset $V_{0}\subset V$ 
such that $V$ is the unique member of ${\cal N}_{0}$ 
which intersects $V_{0}$. 
\end{lemma}
{\bf Proof of Lemma 5.4}. 
Supppose the contrary. 
Let $V$ be a very general member of ${\cal N}_{0}$. 
Let $\eta$ denote the generic point of $V$. 
We define the closed subset $V_{1}$ of $X$ by
\[
V_{1} = \mbox{the closure of}\,\, p(\varphi^{-1}(\varphi (p^{-1}(\eta )))).
\]
By the assumption we see that $\dim V_{1} > \dim V$ 
holds (if there are only finitely many members of ${\cal N}_{0}$ 
which intersect $V$, for a suitable choice of $V_{0}$ 
the assertion is cleary satisfied). 
We note that $V_{1}$ may be reducible. 
Let $S_{1}$ be the closed subset of the closure of $\varphi^{-1}(\varphi (p^{-1}(\eta )))$ defined by 
\[
S_{1}: = \mbox{the closure of}\,\, p^{-1}(\eta ). 
\]
We note that $p^{*}(L,h)$ is numerically trivial on $S_{1}$
by Theorem 4.2, since $(L,h)$ is numerically trivial on 
$V$ and $p(S_{1}) = V$ holds. 
By Lemma 5.1, we see that 
 $p^{*}(L,h)$ is numerically trivial on 
\[
V^{(1)} := \mbox{the closure of $\varphi^{-1}(\varphi (p^{-1}(\eta ))$}, 
\]
since $p^{*}(L,h)$ is numerically trivial on $S_{1}$  and every very general fiber of $\varphi : {\cal V}\longrightarrow {\cal N}^{0}$ by the definition of 
${\cal N}_{0}$.
Here we have used the fact that  the numerical triviality is invariant under modifications, hence we may use the notion of numerical triviality on singular varieties
by Theorem 4.2 (cf. Definition 4.2). 
Again by Theorem 4.2, we  see that $(L,h)$ is numerically trivial on $V_{1}$. 
Since $\dim V_{1}> \dim V$ holds, this contradicts the 
definition of $\nu$. 
This completes the proof of Lemma 5.4. 
{\bf Q.E.D.} \vspace{10mm} \\
Let us continue the proof of Theorem 1.1.
Let ${\cal N}^{\prime}$ be another subvariety of the Hilbert scheme of 
$X$ whose members dominate $X$ and for a very general point $x$ of $X$ there exists a member $V^{\prime}$ of ${\cal N}^{\prime}$ such that 
\begin{enumerate}
\item $x \in V^{\prime}$,
\item $\dim V^{\prime}= \nu$,
\item $(L,h)$ is numerically trivial on $V^{\prime}$.
\end{enumerate} 
Let 
\[
\varphi^{\prime} : {\cal V}^{\prime} \longrightarrow {\cal N}^{\prime}
\]
be the universal family and let 
\[
p^{\prime} : {\cal V}^{\prime} \longrightarrow X
\]
be the natural morphism. 
Then we set  
\[
V_{1}^{\prime} := p^{\prime}((\varphi^{\prime})^{-1}(\varphi^{\prime}((p^{\prime})^{-1}(V))))
\]  
for a very general member $V$ of ${\cal N}^{0}$. 
Repeating the same argument as above we see that 
$(L,h)$ is numerically trivial on $V_{1}^{\prime}$. 
Hence by the definition of $\nu$, we see that 
${\cal N}^{\prime} = {\cal N}^{0}$ holds. 
 
Hence by Lemma 5.4, we see that for a very general point 
$x\in X$, there exists a unique member $V$ of ${\cal N}(x)$ 
such that 
\begin{enumerate}
\item $[V]\in {\cal N}^{0}$, 
\item $(L,h)$ is numerically trivial on $V$,
\item  $\dim V = \nu$.
\end{enumerate} 
Hence there exists a complement $U_{0}$ of
at most  countably many union of proper 
Zariski closed subsets in $X$ such that 
for every $x\in U_{0}$, 
\[
f(x) = [V] \in {\cal N}^{0}, x\in V
\]
is a well defined morphism. 
Hence $f$ defines a rational fibration 
\[
f : X -\cdots\rightarrow Y
\]
by setting 
\[
Y := {\cal N}^{0}. 
\]
If we replace the second condition on $V^{\prime}$,i.e.
$\dim V^{\prime} = \nu$ by $\dim V^{\prime} > 0$, 
by repeating the same argument, we see that 
$V^{\prime}$ is contained in a member of ${\cal N}^{0}$. 
This implies the 3-rd assertion of Theorem 1.1. 
Lemma 5.1 implies the first assertion of Theorem 1.1.
 
By the construction this is the desired fibration. 
This completes the proof of Theorem 1.1. {\bf Q.E.D.}

\section{An algebraic counterpart of the fibration theorem}
An algebraic counterpart of Theorem 1.1 would be the 
following theorem. 
\begin{theorem}
Let $X$ be a normal projective variety and let $L$ be a 
nef line bundle on $X$. 
Then there exists a unique (up to birational equivalence)
 rational fibration 
\[
f : X -\cdots\rightarrow Y
\]
such that 
\begin{enumerate}
\item $f$ is regular over the generic point of $Y$, 
\item  $L$ is numerically trivial on every fibers of $f$,
\item $\dim Y$ is minimal among such fibrations.
\end{enumerate}
\end{theorem}
The proof of the above theorem is essentially same 
as the proof of Theorem 1.1 and is much easier. 
Hence we omit it. 

We should note that the fibrations given by Corollary 1.1 and 
Theorem 6.1 may not be same in general. 
I do not know how to generalize Theorem 6.1 to the case 
that $L$ is pseudoeffective.

Author's address\\
Hajime Tsuji\\
Department of Mathematics\\
Tokyo Institute of Technology\\
2-12-1 Ohokayama, Megro 152-8551\\
Japan \\
e-mail address: tsuji@math.titech.ac.jp
\end{document}